\documentclass{article}
\usepackage{arxiv}
\usepackage[utf8]{inputenc} 
\usepackage[T1]{fontenc}    
\usepackage{hyperref}       
\usepackage{url}            
\usepackage{booktabs}       
\usepackage{amsfonts}       
\usepackage{nicefrac}       
\usepackage{lipsum}
\usepackage[english]{babel}
\usepackage{amsmath,amssymb,amsthm}
\usepackage{mathtools}
\usepackage{algorithm, algorithmic}
\usepackage{pifont}
\usepackage{cases}
\usepackage{graphicx}
\usepackage{stackengine}    
\usepackage{titlesec}
\usepackage[none]{hyphenat}
\usepackage{subfigure}
\usepackage[labelsep=period]{caption}
\usepackage{hyperref}
\hypersetup{
  colorlinks   = true,    
  urlcolor     = green,   
  linkcolor    = blue,    
  citecolor    = red      
}

\newtheorem{theorem}{Theorem}[section]

\newtheorem{remark}[theorem]{Remark}
\newtheorem{proposition}[theorem]{Proposition}

\def\keywordname{{\bfseries Keywords:}}%
\def\keywords#1{\par\addvspace\medskipamount{\rightskip=0pt plus1cm
\def\and{\ifhmode\unskip\nobreak\fi\ $\cdot$ }\noindent\keywordname\enspace\ignorespaces#1\par}}

\def\subjclassname{{\bfseries Mathematics Subject Classification (2000):}}%
\def\subjclass#1{\par\addvspace\medskipamount{\rightskip=0pt plus1cm
\def\and{\ifhmode\unskip\nobreak\fi\ $\cdot$ }\noindent\subjclassname\enspace\ignorespaces#1\par}}

\newcommand\br[1]{\left(#1\right)}

\newcommand\nbr[1]{\left\|#1\right\|}

\titlelabel{\thetitle. }

\begin{document}

\renewcommand{\figurename}{Fig.}

\title{Stopping Rules for Gradient Methods for Non-Convex Problems with Additive Noise in Gradient}

\author{
	Fedor Stonyakin \\
    Moscow Institute of Physics and Technology \\
	Moscow, Russia \\
	V.\,I.\,Vernadsky Crimean Federal University \\
	Simferopol, Russia \\
	\texttt{fedyor@mail.ru} \\
	\And
	Ilya Kuruzov \\
	Moscow Institute of Physics and Technology \\
	Moscow, Russia \\
	\texttt{kuruzov.ia@phystech.edu}
	\And
	Boris Polyak \\
	Moscow Institute of Physics and Technology \\
	Moscow, Russia \\
	Institute for Control Sciences \\
	Moscow, Russia \\
	\texttt{boris@ipu.ru}
}

\maketitle

\phantomsection\addcontentsline{toc}{section}{Abstract}
\begin{abstract}
    We study the gradient method under the assumption that an additively inexact gradient is available for, generally speaking, non-convex problems. The non-convexity of the objective function, as well as the use of an inexactness specified gradient at iterations, can lead to various problems. For example,  the trajectory of the gradient method may be far enough away from the starting point. On the other hand, the unbounded removal of the trajectory of the gradient method in the presence of noise can lead to the removal of the trajectory of the method from the desired exact solution.  The results of investigating the behavior of the trajectory of the gradient method are obtained under the assumption of the inexactness of the gradient and the condition of gradient dominance. It is well known that such a condition is valid for many important non-convex problems. Moreover, it leads to good complexity guarantees for the gradient method. A rule of early stopping of the gradient method is proposed.  Firstly, it guarantees achieving an acceptable quality of the exit point of the method in terms of the function. Secondly, the stopping rule ensures a fairly moderate distance of this point from the chosen initial position.  In addition to the gradient method with a constant step, its variant with adaptive step size is also investigated in detail, which makes it possible to apply the developed technique in the case of an unknown Lipschitz constant for the gradient. Some computational experiments have been carried out which demonstrate effectiveness of the proposed stopping rule for the investigated gradient methods.

    \keywords{Non-convex Optimization
    \and Polyak-Łojasiewicz Condition
    \and Inexact Gradient
    \and Stopping Rule
    \and Adaptive Method}
    \subjclass{49M37 \and 90C25 \and 65K05}
\end{abstract}

\newpage

\section{Introduction}

Gradient methods are relatively simple, and they require a low iteration cost as well as a small amount of memory, which explains their popularity. In data analysis,  non-convex problems often arise under the standard assumption that the gradient of the objective function $f$ is Lipschitz-continuous with some constant $L > 0 $ (or in other words, the function $f$ is $L$-smooth):
\begin{equation}\label{Lsmooth}
\|\nabla f(x) - \nabla f(y)\| \leqslant L \|x - y\| \quad \forall x, y \in \mathbb{R}^n,
\end{equation}
where $\|\cdot\|$ (here and everywhere in the paper) denotes the Euclidean norm. For these problems, by applying the gradient-type methods, the generated sub-sequence of points, generated by applying the gradient-type methods, converges to the zero value of the $\|\nabla f(x)\|$. For this fact we have the following known result (see, for example, \cite{Gasnikov,Polyak1983}).
\begin{theorem}
Let $f$ be an $L$-smooth function. Let us consider the gradient method
\begin{equation}
\label{eqf_3}
    x_{k+1}=x_{k}-\frac{1}{L}\nabla f(x_{k})
\end{equation}
for the following optimization problem:
\begin{equation}\label{optim_problem}
    \min_{x\in\mathbb{R}^{n}} f(x).
\end{equation}
Then the following inequality holds:
\begin{equation}\label{eqf_4}
\min_{k=0,\ldots,N-1}\|\nabla f(x_{k})\|\leqslant\sqrt{\frac{2L(f(x_{0})-f(x_{*}))}{N}},
\end{equation}
where $x_0$ is a starting point of the method and $x_*$ is one of the exact solutions of the problem~\eqref{optim_problem}.
\end{theorem}

Let $f$ be an $L$-smooth function and its gradient satisfy the Polyak-Łojasiewicz condition (for brevity, we will write PL-condition) for some constant $\mu>0$ \cite{Polyak1963} (see also the recent papers \cite{Karimi,Belkin}, and the references therein):
\begin{equation}\label{f1}
f(x)-f^{*}\leqslant\frac{1}{2\mu}\|\nabla f(x)\|^{2}\quad \forall\,x\in\mathbb{R}^{n},
\end{equation}
where $f^* = f(x_*)$ is the value of the function $f$ at one of the exact solutions $x_{*}$ of the optimization problem under consideration.
Then the Gradient Descent Method converges at the rate of a geometric progression
\begin{equation}\label{1.16}
f(x_N)-f^*\leqslant\left(1-\frac{\mu}{L}\right)^N (f(x_0)-f^*) \leqslant\exp\left(-\frac{\mu}{L}N\right)(f(x_0)-f^*),
\end{equation}
\begin{equation}\label{fthm1}
    \|x_{*}-x_{0}\|\leqslant \frac{\sqrt{2L(f(x_{0})-f^{*})}}{\mu}.
\end{equation}

From \cite{Karimi} it is known that the PL-condition~\eqref{f1} implies the following so-called quadratic growth condition:

$$
f(x) - f^{*} \geqslant \frac{\mu}{2} \inf_{x_*}\|x - x_*\|^2 \quad \forall x \in \mathbb{R}^n,
$$
whence one can obtain that~\eqref{1.16} means that the Gradient Descent method also converges in argument at the rate of a geometric progression
$$
\inf_{x_*}\|x_N - x_*\|^2 \leqslant \frac{2}{\mu}\exp\left(-\frac{\mu}{L}N\right)(f(x_0)-f^*).
$$

It is worth noting that the gradient dominance condition~\eqref{f1} is certainly holds for a strongly convex objective function $f$. However, there are known examples, where PL-condition holds, but one cannot be sure even that $f$ is convex (see, for example, \cite{Polyak2017}). So from \cite{Gasnikov}, we can consider the problem of finding some solution to a system of nonlinear equations $g(x)=0$ (written in a vector form), where $g:\mathbb{R}^n\rightarrow\ \mathbb{R}^m$, $m\leqslant n$ the problem of finding some solution to this system.

Let us introduce the Jacobian matrix $J(x)=\frac{\partial g(x)}{\partial x}=\nbr{\frac{\partial g_i(x)}{\partial x_j}}_{i,j=1 }^{m,n}$ of the mapping $g$ and assume that there exists $\mu>0$ such that for all $x\in\mathbb{R}^n$ the Jacobian matrix is uniformly non-singular, i.e.
$\lambda_{min}\br{J(x)\left[J(x)\right]^\top}\geqslant\mu.$ In this case, the function $f(x)=\nbr{g(x)}^2$ satisfies condition~\eqref{f1} for an arbitrary $x_*$ such that $f(x_*)=0$, i.e. $g( x_*)=0$ \cite{Nesterov2006}.
We would like to mention separately the review \cite{Belkin}, which describes in detail a deep learning-motivated example of a non-linear equation-related minimization problem with over-parametrization for a non-convex smooth function with PL-condition.

\subsection{The formulation of the problem}

In this paper, we consider the problem of minimizing the function $f$ which satisfies PL-condition~\eqref{f1} and has $L$-Lipschitz continuous gradient with some constant $L>0$
\begin{equation}\label{f2}
\|\nabla f(x)-\nabla f(y)\|\leqslant L\|x-y\|\quad \forall\,x,\,y\in\mathbb{R}^{n}.
\end{equation}
 We suppose that the method has access not to the exact, but to the approximate value of the gradient $\widetilde{\nabla}f(x)$ at any requested point $x$, which means the following

\begin{equation}
\label{inexact_grad}
    \nabla f(x)=\widetilde{\nabla}f(x)+ v(x),\quad\text{and}\quad \|v(x)\|\leqslant\Delta
\end{equation}
for some fixed $\Delta>0$. Then~\eqref{f1} means that
\begin{equation}\label{f3}
f(x)-f^{*}\leqslant\frac{1}{\mu}(\|\widetilde{\nabla}f(x)\|^{2}+\Delta^{2})\quad \forall\,x\in\mathbb{R}^{n}.
\end{equation}
Therefore, $\|\widetilde{\nabla}f(x)\|^{2}+\Delta^{2}\geqslant\mu(f(x)-f^{*})$, where
 \begin{equation}\label{f4}
 \|\widetilde{\nabla}f(x)\|^{2}\geqslant\mu(f(x)-f^{*})-\Delta^{2}\quad \forall\,x\in\mathbb{R}^{n}.
\end{equation}

It is worth noting that the issue of studying the influence of gradient errors on the estimates of the convergence rate  of the first-order methods attracted many researchers (see, for example, \cite{Polyak1983,s1,DevolderThesis,Artemont,Vasin2021}). However, we will focus on the distinguished class of non-convex problems.
The non-convexity of the objective function of the problem, as well as the use  of an inexactness of the specified gradient at iterations, can lead to various problems. In particular, in the absence of any early stopping rules, divergence of the gradient method trajectory from the starting point can be  quite a large. It is problematic when the initial point of the method already has some appropriate properties. On the other hand, the unlimited divergence of the trajectory of the Gradient Descent method in the stochastic setting can lead to a larger distance from the desired exact solution. Let us describe some situations of this type.

As a simple example of a non-strongly convex function that satisfies the gradient dominance condition, we consider
\begin{equation}\label{simple_ex_1}
f(x) = \langle Ax, x \rangle,
\end{equation}
where $A = \text{diag}(L, \mu, 0)$ is a $3$-order diagonal matrix with exactly two positive entries $L > \mu > 0$.
If for the problem of minimizing the function~\eqref{simple_ex_1} we assume that there is a gradient error $v(x) = (0, 0, \Delta)$ for $\Delta>0$, then for $x_0 = (0, 0, 0)$, $h_k > 0$ and $x_{k+1} = x_k - h_k \widetilde{\nabla}f(x_k)$, we have $\lim\limits_{k \to \infty} \|x_{k+1}\|_2 = \infty$.

Further, we can consider the Rosenbrock function of two variables $x = (x^{(1)}, x^{(2)})$:
$$
f(x) = 100\left(x^{(2)} - \left(x^{(1)}\right)^2\right)^2 + \left(1 - x^{(1)}\right)^2.
$$
Let our method starts from $x_0 = (1, 1) = x_*$. Then at each step of the gradient method, the error of the gradient $v(x_k)$ is such that $x_k^{(2)} = \left(x_k^{(1)}\right)^2$ and without stopping rule the trajectory can go very far from the exact solution $x_*$. Similarly, the trajectory of the gradient method can be unbounded for the objective function of two variables $f(x) = (x^{(2)} - (x^{(1)})^2)^2$.

The purpose of this paper is to study the estimate of the distance $\|x_N - x_0\|$ for points $x_N$ produced by the Gradient Descent method and to propose an early stopping rule that guarantees some compromise, such as a significant divergence of the trajectory from the chosen starting point of the method. Note that the  early stopping rules in iterative procedures are being actively studied for various types of problems.
Apparently, for the first time, the ideology of early stopping of iterations was proposed in \cite{EmelinKrasnos}. This paper is devoted to a technique for the approximate solution of ill-posed or ill-conditioned problems arising during regularization (in the mentioned work,  the authors considered the problem of solving a linear equation).
In this case, an early stop is aimed at overcoming the problem of the potential accumulation of errors in the regularization of the original problem. The topic of our paper is related to well-known approaches related to the early termination of first-order methods in the case of using inexact information about the gradient at iterations (see \cite{Polyak1983}, Ch. 6, paragraph 1, and also, for example, the recent preprint \cite{Vasin2021}). However, the  results known to us for convex (not strongly convex) problems differ from those obtained in this note. The main difference is that usually either the achievement of the worst level in function is guaranteed (compared with the comment after theorem 2, section 1, chapter 6 of \cite{Polyak1983}) or estimates such as $\|x_N - x_*\| \leqslant \|x_0 - x_*\|$ without examining $\|x_N - x_0\|$. Here $\{x_k\}_{k \in \mathbb{N}}$ is the sequence generated by the method, $ x_*$ is the exact solution of the minimization problem closest to the starting point of the method $x_0$.

In this paper, we obtained  Theorem~\ref{main_thm_2} devoted to the Gradient Descent method with a constant step-size with a sufficiently small value of the inexact gradient. It indicates the level of accuracy with respect to the function that can be guaranteed after the proposed early stopping rule is fulfilled. It is important to note that this result can be applied to any $L$-smooth non-convex problem. Further, using PL-condition, this result is refined in Theorem~\ref{thm-main}, which describes the estimate of a sufficient number of iterations to achieve the desired quality of the output point $\widehat{x}$ by the function $f(\widehat{x }) - f^* = O\left(\frac{\Delta^2}{\mu}\right)$. Moreover, it contains an estimate~\eqref{f12} of the distance from $\widehat{x}$ to the starting point $x_0$.
The obtained results  are compared with the well-known distance estimate \cite{Polyak1963} from the starting point $x_0$ to the nearest exact solution $x_*$ (see remark~\ref{rem-nain}).

However, the method with a constant step-size imposes the need to efficiently estimate the Lipschitz constant of the gradient of the objective function, which can be problematic in practice. Moreover, many real problems lead to functions that have not an Lipschitz continuous gradient, and a condition such as~\eqref{Lsmooth} holds for such functions only locally on some subset $Q \subset \mathbb{R}^n$. Therefore, we propose variations of Theorems~\ref{main_thm_2} and~\ref{thm-main}  for the Gradient Descent method with an adaptively selected step-size. This makes it possible to apply an analog of Theorem~\ref{main_thm_2}  with the early stopping rule~\eqref{stop_cond_adaptGD} to an arbitrary non-convex problem without additional conditions. If PL-condition is guaranteed, then the execution of~\eqref{stop_cond_adaptGD} automatically guarantees the achievement of an acceptable quality level for the solution of the problem of minimizing $f$ by the function.

The last section of the paper is devoted to  numerical experiments which explain the purpose of using stopping rule~\eqref{f9} for some specific examples of object functions in problems: logistic regression, Rosenbrock and Nesterov-Skokov functions, quadratic function.

\section{The Proposed Approach and Main Theoretical Results}

\subsection{Variant of the Gradient Descent method with a constant step-size}
\label{sec_const_step}

We assume that the values of the parameters $L>0$ and $\Delta>0$ are known. Also, the Gradient Descent  method of the following form
\begin{equation}\label{f5}
x_{k+1}=x_{k}-\frac{1}{L}\widetilde{\nabla}f(x_{k})
\end{equation}
can be applied to solve the minimization problem of the function $f$. In view of~\eqref{f2} for the method~\eqref{f5}, we get
\begin{equation*}
    \begin{aligned}
    f(x_{k+1}) & \leqslant f(x_{k})+\langle \nabla f(x_{k}), x_{k+1}-x_{k}\rangle+\frac{L}{2}\|x_{k+1}-x_{k}\|^{2}  \\&
    = f(x_{k}) - \frac{1}{L}\langle\nabla f(x_{k}), \widetilde{\nabla} f(x_{k}) \rangle+ \frac{1}{2L} \|\widetilde{\nabla} f(x_{k})\|^{2} \\&
    = f(x_{k}) + \frac{1}{2L}\left(\|\nabla f(x_{k})\|^{2} - 2 \langle \nabla f(x_{k}), \widetilde{\nabla} f(x_{k}) \rangle+ \|\widetilde{\nabla} f(x_{k})\|^{2}\right) - \frac{\|\nabla f(x_{k})\|^{2}}{2L}  \\&
    = f(x_{k}) + \frac{1}{2L} \|\nabla f(x_{k}) - \widetilde{\nabla} f(x_{k})\|^2 - \frac{\|\nabla f(x_{k})\|^{2}}{2L} \\&
    \leqslant f(x_{k}) + \frac{\Delta^2}{2L} -\frac{1}{2L}\|\nabla f(x_{k})\|^{2},
    \end{aligned}
\end{equation*}
i.e.
\begin{equation}\label{f6}
f(x_{k+1}) - f(x_{k})\leqslant \frac{\Delta^{2}}{2L} - \frac{1}{2L}\|\nabla f(x_{k})\|^{2}.
\end{equation}
Summing up inequalities~\eqref{f6} over $k = \overline{0, N-1}$ leads us to an estimate
\begin{equation}\label{eqnormgrad}
\min\limits_{k = 0,\ldots, N-1}\|\nabla f(x_{k})\| \leqslant \sqrt{\Delta^2 + \frac{2L(f(x_{0})-f(x_{*}))}{N}} \leqslant \Delta + \sqrt{\frac{2L(f(x_{0})-f(x_{*}))}{N}}.
\end{equation}
Note that, in contrast to~\eqref{eqf_4}, the estimate~\eqref{eqnormgrad} points to the potential divergence of the Gradient Descent method in the case of an additively inexact gradient. Specific examples of such situations were described above.

Taking into account~\eqref{f1} , we get
$$
f(x_{k+1}) - f(x_{k})\leqslant \frac{\Delta^{2}}{2L}-\frac{2\mu(f(x_{k})-f^{*})}{2L} = - \frac{\mu}{L}(f(x_{k})-f^{*})+\frac{\Delta^{2}}{2L},
$$
thus
\begin{align*}
    f(x_{k+1})-f^{*}&\leqslant \left(1- \frac{\mu}{L}\right)(f(x_{k})-f^{*})+\frac{\Delta^{2}}{2L}\\&
    \leqslant\left(1- \frac{\mu}{L}\right)^{k+1}(f(x_{0})-f^{*})+\frac{\Delta^{2}}{2L}\left(1+1-\frac{\mu}{L}+\cdots +\left(1-\frac{\mu}{L}\right)^{k}\right) \\&
    <\left(1- \frac{\mu}{L}\right)^{k+1}(f(x_{0})-f^{*})+\frac{\Delta^{2}}{2\mu},
\end{align*}
i.e.
\begin{equation}\label{f7}
f(x_{k+1})-f^{*}\leqslant \left(1- \frac{\mu}{L}\right)^{k+1}(f(x_{0})-f^{*})+\frac{\Delta^{2}}{2\mu}.
\end{equation}

\begin{remark}
It is important to note that bounds~\eqref{eqnormgrad} and~\eqref{f7} cannot be improved for the Gradient Descent method with an additively inexact gradient in the general case. For example, the lower estimates of accuracy with respect to the function $O\left(\frac{\Delta^2}{2\mu}\right)$  are known even on the class of strongly convex functions (see, for example, section 2.11.1 of the manual \cite{Vorontsova2021}, as well as references therein). In this regard, we consider the following example:
\begin{equation}\label{f_1_1}
\min_{x\in\mathbb{R}^{n}}f(x):=\frac{1}{2}\sum_{i=1}^{n}\lambda_{i}\left(x^{(i)}\right)^{2},
\end{equation}
where $0\leqslant\mu=\lambda_{1}\leqslant\lambda_{2}\leqslant\ldots\leqslant\lambda_{n}=L$, such that $L\geqslant 2\mu$. The exact solution of problem~\eqref{f_1_1}  is $x_{*}=0 \in \mathbb{R}^n$. Suppose that an inexact gradient is available at the current point of the feasible area. Besides, the error is only in the calculation of the first component of the gradient. That is, instead of $\partial f(x)/\partial x^{(1)}=\mu x^{(1)}$, we have only $\widetilde{\partial} f(x)/\partial x ^{(1)}=\mu x^{(1)}-\Delta$, for some $\Delta >0$. Then for the simplest Gradient Descent method~\eqref{f5}
one can obtain that for $x^{(1)}_{0}\geqslant0$ and sufficiently large $k\in\mathbb{N}$  $(k\gg L/\mu)$ the following inequality holds:
\begin{equation}\label{f_1_2}
x^{(1)}_{k}\geqslant\frac{\Delta}{L}\frac{1-(1-\mu/L)^{k}}{1-(1-\mu/L)}\simeq\frac{\Delta}{\mu}.
\end{equation}
Therefore, $
f(x_{k})-f(x_{*})\gtrsim\frac{\Delta^{2}}{2\mu}.
$
\end{remark}
Further, in view of
$$
\|\nabla f(x_{k})\|^2 \geqslant  \frac{\|\widetilde{\nabla} f(x_{k})\|^2}{2} - \Delta^2,
$$
from~\eqref{f6}, the inexact gradient satisfies the following inequality:
$$
f(x_{k+1}) - f(x_{k})\leqslant \frac{\Delta^{2}}{2L} - \frac{1}{2L}\left(\frac{\|\widetilde{\nabla} f(x_{k})\|^2}{2} - \Delta^2\right),
$$
whence we have
\begin{equation}\label{f61}
f(x_{k+1}) - f(x_{k})\leqslant \frac{\Delta^{2}}{L} - \frac{1}{4L}\|\widetilde{\nabla} f(x_{k})\|^{2}.
\end{equation}

Inequality~\eqref{f61} shows that if the value $\|\widetilde{\nabla} f(x_{k})\|$ is sufficiently large, it can be guaranteed that $f(x_{k+1}) < f( x_{k})$. Thus, for any $C > 2$, an alternative arises: either the inequality $\|\widetilde{\nabla} f(x_{k})\| \leqslant C\Delta$ holds, or
$$f(x_{k+1}) - f(x_{ k}) < -\frac{\Delta^2}{L}\left(\frac{C^2}{4} - 1\right).
$$
In the first case, the inequality $\|\widetilde{\nabla} f(x_{k})\| \leqslant C\Delta$ guarantees the achievement of an acceptable quality of the output point $x_k$ with respect to the function due to PL-condition. In the second case, we can guarantee the decreasing with respect to the function for $C>2$.

So, it is possible to get $x_k$ such that the value of $f(x_k)$ is close enough to the minimum $f^*$. For definiteness, let us choose $C = \sqrt{6}$ (to get a “convenient” coefficient) and consider 2 scenarios:
\begin{enumerate}
\item $\|\widetilde{\nabla} f(x_{k})\|> \Delta \sqrt{6}$, whence, taking~\eqref{f61} into account, we obtain the inequality
\begin{equation}\label{f8}
f(x_{k+1}) - f(x_{k})<-\frac{\Delta^{2}}{2L}.
\end{equation}
\item
\begin{equation}\label{f9}
\|\widetilde{\nabla} f(x_{k})\|\leqslant \Delta\sqrt{6},
\end{equation}
\end{enumerate}
whence, in view of~\eqref{f3} we have
\begin{equation}\label{f10}
f(x_{k})-f^{*}\leqslant\frac{7\Delta^{2}}{\mu}.
\end{equation}

Let us consider estimate~\eqref{f10} acceptable for the function level and {\it agree to terminate process~\eqref{f5} if~\eqref{f9} is satisfied.}

Let us investigate an alternative situation in which for any $k=0,1,\ldots,N-1$, it is true that  $\|\widetilde{\nabla} f(x_{k})\|>\Delta \sqrt{6}$ and ~\eqref{f8} holds, where
$$
f(x_{0})-f(x_{N})=\sum_{k=0}^{N-1}(f(x_{k})-f(x_{k+1}))>\frac{N\Delta^{2}}{2L},
$$
i.e.
$
N < \frac{2L}{\Delta^{2}}(f(x_{0})-f^{*}),
$
which indicates the end of the process. Thus, we have the following result.

\begin{theorem}\label{main_thm_2}
Let stopping criterion~\eqref{f9} be satisfied for the first time at the $N$-th iteration of the Gradient Descent method~\eqref{f5}. Then the output point $\widehat{x} = x_{N}$ is guaranteed to satisfy the inequality
$$
f(\widehat{x})-f^{*}\leqslant\frac{7\Delta^{2}}{\mu}.
$$
In this case, the following estimate for the number of iterations before stopping criterion is valid
\begin{equation}\label{finit_proc}
N < \frac{2L}{\Delta^{2}}(f(x_{0})-f^{*}).
\end{equation}
\end{theorem}

It is clear that for a small value of the parameter $\Delta > 0$, the right-hand side of inequality~\eqref{finit_proc} leads to a significantly overestimated number of iterations. At the same time, the conducted computational experiments  (see Section 4, below) showed no increase in the number of iterations with a significant decrease in $\Delta > 0$ due to the proposed early stopping rule~\eqref{f9}.

However, in the case of a known $\mu$, the estimate for the number of steps $N$ can be improved if the quality in~\eqref{f10} is assumed to be sufficient. Using inequality~\eqref{f61} we get
$
\frac{1}{4L}\|\widetilde{\nabla} f(x_{k})\|^{2}\leqslant\frac{\Delta^{2}}{L} +f(x_{k})-f(x_{k+1}),
$ and due to $\widetilde{\nabla} f(x_{k})=L(x_{k}-x_{k+1})$, we have the following estimation for every $k\geqslant 0$:
\begin{equation*}
    \begin{aligned}
      \|x_{k+1}-x_{k}\|^{2} & \leqslant\frac{4\Delta^{2}}{L^{2}} +\frac{4(f(x_{k})-f(x_{k+1}))}{L} \\&
      \leqslant \frac{4\Delta^{2}}{L^{2}} +\frac{4(f(x_{k})-f^{*})}{L}
      \\&
      \leqslant \frac{4\Delta^{2}}{L^{2}}+\frac{4\Delta^{2}}{\mu L}+\frac{4}{L}\left(1-\frac{\mu}{L}\right)^{k}(f(x_{0})-f^{*}).
    \end{aligned}
\end{equation*}
Whence one can obtain the final estimation:
 $$
\|x_{k+1}-x_{k}\|  \leqslant 2\Delta\sqrt{\frac{1}{L^{2}}+\frac{1}{\mu L}}+ 2 \left(1-\frac{\mu}{L}\right)^{\frac{k}{2}}\sqrt{\frac{f(x_{0})-f^{*}}{L}}.
 $$

Next, summing the inequalities above for $k=0\dots N-1$, we have
\begin{equation}\label{f11}
\|x_{N}-x_{0}\|\leqslant\sum_{k=0}^{N-1}\|x_{k+1}-x_{k}\|  \leqslant 2N\Delta\sqrt{\frac{1}{L^{2}}+\frac{1}{\mu L}}+ 2 \sum_{k=0}^{N-1}\left(1-\frac{\mu}{L}\right)^{\frac{k}{2}}\sqrt{\frac{f(x_{0})-f^{*}}{L}}.
\end{equation}

If at some step~\eqref{f9} is satisfied, then the required accuracy by function~\eqref{f10} will be achieved. Therefore, we estimate $N$ in an alternative situation (\eqref{f9} does not hold for all $k=0,1,\ldots,N-1$). We use  inequality~\eqref{f7} and impose the requirement that the level of approximation with respect to the function $f(x_{N})-f^{*}\leqslant\frac{7\Delta^{2}}{\mu}.$ In view of~\eqref{f7}, it suffices to require that
$$
\left(1-\frac{\mu}{L}\right)^{N}(f(x_{0})-f^{*})\leqslant\frac{6\Delta^{2}}{\mu},
$$
or
$$
\left(1-\frac{\mu}{L}\right)^{N}\leqslant e^{-\frac{\mu N}{L}}\leqslant\frac{6\Delta^{2}}{\mu(f(x_{0})-f^{*})},
$$
where
$
N\leqslant\left\lceil \frac{L}{\mu}\ln\frac{\mu(f(x_{0})-f^{*})}{6\Delta^{2}}\right\rceil.
$
In this case,~\eqref{f11} takes the following form:
$$
\|x_{N}-x_{0}\|\leqslant\frac{2\Delta}{\mu}\sqrt{1+\frac{L}{\mu}}\left\lceil\ln\frac{\mu(f(x_{0})-f^{*})}{6\Delta^{2}}\right\rceil+\frac{4\sqrt{L(f(x_{0})-f^{*})}}{\mu}.
$$

\begin{theorem}\label{thm-main}
Let one of the following alternatives hold:
\begin{enumerate}
    \item The Gradient Descent method~\eqref{f5} works $N_*$ steps where $N_*$ is such that
\begin{equation}\label{equat_estimN*}
N_{*}=\left\lceil \frac{L}{\mu}\ln\frac{\mu(f(x_{0})-f^{*})}{6\Delta^2} \right\rceil.
\end{equation}
\item For some $N \leqslant N_*$, at the $N$-th iteration of the method~\eqref{f5}, stopping criterion~\eqref{f9} is satisfied for the first time.
\end{enumerate}
Then for the output point $\widehat{x}$ ($\widehat{x} = x_{N}$ or $\widehat{x} = x_{N_*}$) of the method~\eqref{f5}, the following inequalities hold:
$$
f(\widehat{x})-f^{*}\leqslant\frac{7\Delta^{2}}{\mu},
$$
\begin{equation}\label{f12}
\|\widehat{x}-x_{0}\|\leqslant\frac{2\Delta}{\mu}\sqrt{1+\frac{L}{\mu}}\left\lceil\ln\frac{\mu(f(x_{0})-f^{*})}{6\Delta^{2}}\right\rceil+\frac{4\sqrt{L(f(x_{0})-f^{*})}}{\mu}.
\end{equation}
\end{theorem}

\begin{remark}\label{rem-notmu}
Since it is often difficult to estimate the value of the parameter $\mu$ and usually  $f^*$ is not known, the estimate of the number of iterations~\eqref{equat_estimN*} is difficult to use in practice. If the implementation works only according to stopping rule~\eqref{f9}, then we can only confirm an upper bound on the number of iterations of the form~\eqref{finit_proc}, but in this case we cannot guarantee~\eqref{f12}. However, the estimate~\eqref{f11}  remains relevant. Moreover, the estimation of the value $\|\widehat{x}-x_{0}\|$ can be refined if the value of the parameter $\mu$ is not available. Indeed, in view of~\eqref{f61} for the Gradient Descent method~\eqref{f5} with a constant step-size it holds that
$\frac{1}{4L}\|\widetilde{\nabla} f(x_{k})\|^{2} \leqslant  \frac{\Delta^{2}}{L} + f(x_{k}) - f(x_{k+1}).$
Whence we have
$\|x_{k+1}-x_{k}\|^{2} \leqslant  \frac{4\Delta^{2}}{L^2} + \frac{4(f(x_{k}) - f(x_{k+1})}{L},$
i.e. $
\|x_{k+1}-x_{k}\| \leqslant  \frac{2\Delta}{L} + 2\sqrt{\frac{(f(x_{k}) - f(x_{k+1})}{L}}.
$
Further, after summing the inequalities above over $k = \overline{0, N-1}$, we have:
\begin{equation*}
    \begin{aligned}
    \|x_{0}-x_{N}\| & \leqslant\sum^{N-1}_{k=0}\|x_{k}-x_{k+1}\|\leqslant\frac{2N\Delta}{L} + 2\sum^{N-1}_{k=0}\sqrt{\frac{f(x_{k})-f(x_{k+1})}{L}} \\&
    \leqslant \frac{2N\Delta}{L}+ \sqrt{N}\sqrt{\sum^{N-1}_{k=0}\frac{f(x_{k})-f(x_{k+1})}{L}} \\&
    = \frac{2N\Delta}{L}+ 2\sqrt{N}\sqrt{\frac{f(x_{0})-f(x_{N})}{L}} \\&
    \leqslant \frac{2N\Delta}{L}+ 2\sqrt{N}\sqrt{\frac{f(x_{0})-f^*}{L}}.
    \end{aligned}
\end{equation*}
It is clear that for small values of the error $\Delta>0$ the following inequality
$$
\|x_{0}-x_{N}\|\leqslant  \frac{2N\Delta}{L}+ 2\sqrt{N}\sqrt{\frac{f(x_{0})-f^*}{L}}
$$
may turn out to be worse than~\eqref{f12}. Taking into account~\eqref{finit_proc}, we get 
\begin{equation*}
    \begin{aligned}
    \|x_{0}-x_{N}\| & \leqslant\frac{2\Delta}{L} \cdot\frac{2L(f(x_{0})-f^{*})}{\Delta^{2}} + 2\sqrt{\frac{2L(f(x_{0})-f^{*})}{\Delta^{2}}\cdot\frac{f(x_{0})-f^{*}}{L}}
    \\& = \frac{4 + 2\sqrt{2}}{\Delta}(f(x_{0})-f^{*}).
    \end{aligned}
\end{equation*}
\end{remark}

\begin{remark}\label{rem-nain}
By~\eqref{fthm1} and~\eqref{f12} the quantity $\|\widehat{x}-x_{0}\|$ can be comparable with $\|x_{*}-x_{0}\|$ for a sufficiently small $\Delta>0$.
\end{remark}

\begin{remark}
In view of~\eqref{f12}, it suffices to require that conditions~\eqref{f1} and~\eqref{f2} are satisfied only in $R$--neighborhood of the $x_{0}$, where
$$
R=\frac{2\Delta}{\mu}\sqrt{1+\frac{L}{\mu}}\left\lceil\ln\frac{\mu(f(x_{0})-f^{*})}{6\Delta^{2}}\right\rceil+\frac{4\sqrt{L(f(x_{0})-f^{*})}}{\mu}.
$$
\end{remark}

\subsection{Some Variant of the Gradient Descent Method with an Adaptive Step-Size Policy}

In many applied optimization problems, it is difficult to estimate the Lipschitz constant of the gradient of the objective function. For example, the well-known Rosenbrock function and its multidimensional generalizations (for example, the Nesterov-Skokov function \cite{NestSkok}) have only a locally Lipschitz gradient. Thus, it is impossible to estimate for them the Lipschitz constant of the gradient without additional restrictions on the domain in which the method operates. Therefore, we present a generalization of the universal gradient method from  \cite{NesterovUniversal} for working with an inexact gradient of the  functions satisfying PL-condition.

For $L$-smooth functions we have the following well-known inequality:
$$f(x)\leqslant f(y)+\langle\nabla f(y), x-y\rangle + \frac{L}{2}\|x-y\|^2, \quad \forall x,y\in\mathbb{R}^n.
$$
For the inexact gradient~\eqref{inexact_grad} we can get a similar inequality:
\begin{align*}
    f(x)\leqslant f(y) + \langle\widetilde{\nabla} f(y), x-y\rangle + L\|x-y\|^2 +\frac{\Delta^2}{2L}, \quad  \forall x,y\in\mathbb{R}^n.
\end{align*}

This inequality contains an exact calculation of the value of the function $f$ at an arbitrary point from the $\text{dom} f$. For most important applications with an inexact gradient, we do not have an opportunity to make such a calculation. An important example of such problems is some optimization problems in the Hilbert space \cite{Vasilyev} and, in a particular case, inverse problems \cite{Kabanikhin}. Therefore, further, we will discuss the possibility of using an inexact function value when checking the iteration exit criterion.

Let us assume that we can calculate the inexact value $\widetilde{f}$ of the function $f$ at any point $x$, so that
\begin{equation}
    \label{inexactf_cond}
    |f(x)-\widetilde{f}(x)|\leqslant \delta.
\end{equation}

Then we have the following inequality:
\begin{equation}
    \label{smooth_cond_Lk}
    \widetilde{f}(x)\leqslant \widetilde{f}(y)+\langle\widetilde{\nabla} f(y), x-y\rangle + L\|x-y\|^2 + \frac{\Delta ^2}{2L} + 2\delta,   \;\;\; \forall x,y\in\mathbb{R}^n.
\end{equation}

Further, when $\mu$ is known, we select the constant $L$ in such a way that~\eqref{smooth_cond_Lk} is satisfied for the points from the neighboring iterations (see Algorithm~\ref{adapt_gd}).

\begin{algorithm}
\caption{Adaptive Gradient Descent with Inexact Gradient.}
\label{adapt_gd}
  \begin{algorithmic}[1]
  \REQUIRE $L_{\min}\geqslant 0, L_0\geqslant L_{\min}, \delta\geqslant 0, \Delta\geqslant 0$.
  \STATE Set $k:=0$
  \STATE Calculate 
  \begin{equation}
  \label{x_k}
  x_{k+1} = x_k - \frac{1}{2L_k}\widetilde{\nabla}f(x_k)
  \end{equation}
    \STATE If the following inequality holds:
      \begin{equation}
      \label{cond_it}
    \widetilde{f}(x_{k+1})\leqslant \widetilde{f}(x_k)+\langle\widetilde{\nabla} f(x_k), x_{k+1}-x_k\rangle + L_k\|x_{k+1}-x_k\|^2 + \frac{\Delta^2}{2L_k} +2\delta, \end{equation}
    then $k:=k+1$, $L_k := \max\left(\frac{L_{k-1}}{2}, L_{\min}\right)$ and go to Step 2. Otherwise, $L_k:=2L_k$ and go to Step 3.
\RETURN $x_k$
\end{algorithmic}
\end{algorithm}

Similarly to the approach of the method with a constant step-size proposed above, in the case of a sufficiently small inexact gradient
\begin{equation}
\label{stop_cond_adaptGD}
    \|\widetilde{\nabla}f(x_k)\|\leqslant 2\Delta
\end{equation}
we agree to interrupt Algorithm~\ref{adapt_gd}. In this case, according to~\eqref{f3} we can guarantee that $f(x_k)-f^*\leqslant \frac{5\Delta^2}{\mu}$.

An alternative case, where condition~\eqref{stop_cond_adaptGD}  is not satisfied, can be investigated similarly to the constant step-size case in Section~\ref{sec_const_step}. A detailed proof is given in Appendix~\ref{app:adaptL_inexactf_proof}.

The theoretical results about the operation of Algorithm~\ref{adapt_gd}  are presented in the following theorem.
\begin{theorem}
\label{theorem:adaptL_inexactf}
Suppose $f(x)$ satisfies PL-condition~\eqref{f1} and conditions~\eqref{inexactf_cond}, $\Delta^2 \geqslant 16L\delta$ hold. Let  the parameter $L_{\min}$ in Algorithm~\ref{adapt_gd} be such that $L_{\min}\geqslant\frac{\mu}{4}$ and one of the following alternatives holds:
\begin{enumerate}
    \item Algorithm~\ref{adapt_gd} works $N_*$ steps where $N_*$ is such that
\begin{equation}
    \label{N_star_adapt_gd}
    N_*= \left\lceil\frac{8L}{\mu}\log \frac{\mu (f(x_0)-f^*)}{\Delta^2}\right\rceil.
\end{equation}
\item For some $N \leqslant N_*$, at the $N$-th iteration of Algorithm~\ref{adapt_gd},  stopping criterion~\eqref{stop_cond_adaptGD} is satisfied for the first time.
\end{enumerate}
Then for the output point $\widehat{x}$ ($\widehat{x} = x_{N}$ or $\widehat{x} = x_{N_*}$) of Algorithm~\ref{adapt_gd}, we have the following inequalities
$$
f(\widehat{x})-f^{*}\leqslant\frac{5\Delta^{2}}{\mu},
$$
\begin{equation}
\label{dist}
    \|\widehat{x}-x_0\|\leqslant 8\frac{\Delta}{\mu} \sqrt{\frac{1}{2}\gamma^2  + 4\gamma\frac{L}{\mu }} \log \frac{\mu (f(x_0)-f^*)}{\Delta^2} + 16\frac{\sqrt{\gamma L(f(x_0)-f^*)}}{\mu},
\end{equation}
where $\gamma=\frac{L}{L_{\min}}$. Also, the total number of calls to the subroutine for calculating inexact values of the objective function and step~\eqref{x_k} is not more than $2N+\log\frac{2L}{L_0}.$
\end{theorem}

As we can see, estimate~\eqref{dist} from Theorem~\ref{theorem:adaptL_inexactf} for the Gradient Descent with an adaptive step-size differs significantly from  estimate~\eqref{f12} from Theorem~\ref{thm-main} for the method with a constant step-size, namely, by  the presence of the factor $ \gamma$. In the worst case, the ratio of these two estimates can be $O\left(\frac{L}{\mu}\right)$. However, as it will be shown in experiments, the distances $\|\widehat{x} - x_0\|$ for the methods differ insignificantly. In addition, note, that Algorithm~\ref{adapt_gd} uses subroutines for finding the inexact value of the objective function more often than the gradient method with a constant step. But the number of calls to these subroutines in adaptive Algorithm~\ref{adapt_gd} is not more than  $2N+\log\frac{2L}{L_0}$. This means that the "cost" of an iteration of the adaptive algorithm is on average comparable to about two iterations of the non-adaptive method~\eqref{f5}. At the same time, the accuracy achieved by the proposed methods is also approximately equal.

\begin{remark}
Note that  condition~\eqref{stop_cond_adaptGD} is satisfied for any $L_k\geqslant L$. By construction, we obtain that $L_k\leqslant 2L$. In the estimates above, the quantity $2L$ estimates the maximum value of the parameter $L_k$. The estimates above remain valid if $L$ is replaced by $\frac{1}{2}\max_{j\leqslant k}L_j$ and $\gamma$ by $\frac{\max_{j\leqslant k}L_j}{2\min_{j\leqslant k} L_j}.$ Similarly, we can replace the algorithm parameter $L_{\min}$ with $\min_{j\leqslant k}L_j$.
\end{remark}

\begin{remark}
Note that the estimate for the number of iterations~\eqref{N_star_adapt_gd} in Theorem~\ref{theorem:adaptL_inexactf} indicates the finiteness of the process, but it is strongly overestimated. In practice, the following relation is a more interesting:
\begin{equation*}
    N_*= \left\lceil\frac{4\widehat{L}}{\mu}\log \frac{\mu (f(x_0)-f^*)}{\Delta^2}\right\rceil,
\end{equation*}
where $\widehat{L}=\frac{\mu}{4}\frac{1}{1-\left(\prod_{j=0}^{N_*-1} \left(1-\frac{ \mu}{4L_j}\right)\right)^{\frac{1}{N_*}}}$ is a parameter depending on the fitted parameters $L_j$  in Algorithm~\ref{adapt_gd}.
 \end{remark}

\begin{remark} \label{remark7}
Also note that we can relax the requirement $L_{\min}\geqslant \frac{\mu}{4}$ to $L_{\min}>0$. In this case, the estimate for the distance from the starting point to the point $x_N$  at the $N$-th iteration (see the proof of the expression~\eqref{distN_adaptL}) has the following form
\begin{equation*}
    \|x_N-x_0\|\leqslant N\Delta \sqrt{\frac{1}{2L_{\min}^2} + \frac{4}{\mu L_{\min}}} + 16\sqrt{\frac{L}{L_{\min}}}\frac{\sqrt{L(f(x_0)-f^*)}}{\mu}.
\end{equation*}
But we can no longer use estimate~\eqref{dist} from Theorem~\ref{theorem:adaptL_inexactf}. In this case, it is possible to evaluate the sufficient number of iterations of Algorithm~\ref{adapt_gd}, assuming that the stopping condition $\|\widetilde{\nabla}f(x_k)\|\leqslant 2\Delta$ is not satisfied. Further, we obtain an estimate (see the proof of~\eqref{N_delta}) for $\Delta^{2}>16L\delta$:
$$
N < \frac{2L}{\Delta^{2}-16L\delta}(f(x_{0})-f^{*}).
$$
\end{remark}

For experimental comparison of the Gradient Descent methods with constant and adaptive steps, a single stopping criterion must be chosen. If we consider criterion~\eqref{f9} for the adaptive Algorithm~\ref{adapt_gd} instead of~\eqref{stop_cond_adaptGD}, then the results of Theorem~\ref{theorem:adaptL_inexactf}  about the number of iterations~\eqref{N_star_adapt_gd} and the estimate of the distance from $x_0$ to $ \widehat{x}$ will remain valid. Thus, criterion~\eqref{f9}  makes it possible to achieve the same theoretical guarantees. Therefore, further in the experimental comparison of the variants of the gradient method  (Algorithms~\ref{adapt_gd} and~\eqref{f5}), we will use stopping criterion~\eqref{f9}.

\section{Numerical Experiments}

\subsection{The quadratic form minimization problem}
\label{sec:exp_qp}

In this section, we compare the number of iterations of   method~\eqref{f5}, required to stop according to  criterion~\eqref{f9},  and the estimate for the number of iterations~\eqref{equat_estimN*} to achieve  estimate~\eqref{f10}. To obtain the theoretical estimate for the number of iterations~\eqref{equat_estimN*} we need the values of the constants $L$ and $\mu$. Therefore, as the first example, we consider a quadratic function for which these constants are easy to calculate.

 As an inexact gradient, we use an exact gradient with random noise~\eqref{inexact_grad}. In our experiments, we consider the following types of the additive inexactness $v(x)$ in~\eqref{inexact_grad}:
\begin{itemize}
    \item \textbf{Random.}  Randomly generated from a uniform distribution, i.e.  $v(x)\sim\mathcal{U}\left(S_1^n(0)\right)$, where $S_1^n(0) $ is the $n$ dimensional sphere with radius 1 at the center 0.
    \item \textbf{Antigradient.}
    $v(x)=-\frac{ \nabla f(x)}{\| \nabla f(x)\|}$.
    \item \textbf{Constant.} $v(x)=v\in\mathbb{R}^n,$ such that $ \|v\|=1$.
\end{itemize}

Let us start with a simple example that allows us to estimate the parameters $L$ and $\mu$.  As shown in \cite{Polyak1963}, the function $f(x)=\frac{1}{2}\langle x, A x\rangle$ satisfies PL-condition if the operator $A$ is non-negative definite and its spectrum is separated from zero. In such a case, $\mu$ is the smallest nonzero eigenvalue of the matrix $A$. At the same time, the Lipschitz constant of the gradient is the largest eigenvalue of the matrix $A$. Thus, we consider the following problem of quadratic programming:
\begin{equation}\label{Q_function}
    \min_{x\in\mathbb{R}^n} \sum\limits_{j=k+1}^n d_j x_j^2,
\end{equation}
where $k$ is the number of zero eigenvalues of the matrix $A$, and $d_j$ are some positive constants. Thus, we have a quadratic form with a non-negative definite diagonal matrix. In this case, we can explicitly find the constants $\mu=\min\limits_{j=\overline{k+1,n}}(d_j), L=\max\limits_{j=\overline{k+1 ,n}}(d_j)$.

In the conducted experiments, we take $L=1$ and change $\mu$ from $0$ to $1$. The parameters $d_j$ will be taken uniformly random from the interval $[\mu, L]$.
We take the dimension $n=100$ and $k=10$ of zero eigenvalues. Let us compare the required number of iterations  to achieve condition~\eqref{f9} and the estimate of $N_*$ from Theorem~\ref{thm-main}. As an inexactness, we will take \textbf{Random} noise $v(x)$. The results for problem~\eqref{Q_function} are presented in Table~\ref{tab:k_N}.

\begin{table}[ht]
    \centering
    \begin{tabular}{|c|c|c|c||c|c|c|c|}
      \hline
      $\mu$ & $\Delta$ & $N$ & $N_*$&$\mu$ & $\Delta$ & $N$ & $N_*$ \\
      \hline
0.01 & \begin{tabular}{@{}c@{}} $10^{-7}$ \\ $10^{-4}$ \\ $10^{-1}$ \end{tabular}&\begin{tabular}{@{}c@{}} $1528$ \\ $841$ \\ $155$ \end{tabular}&\begin{tabular}{@{}c@{}} $3817$ \\ $2436$ \\ $1054$ \end{tabular}& 0.1 & \begin{tabular}{@{}c@{}} $10^{-7}$ \\ $10^{-4}$ \\ $10^{-1}$ \end{tabular}&\begin{tabular}{@{}c@{}} $169$ \\ $104$ \\ $40$ \end{tabular}&\begin{tabular}{@{}c@{}} $406$ \\ $267$ \\ $129$ \end{tabular}\\
\hline
0.9 & \begin{tabular}{@{}c@{}} $10^{-7}$ \\ $10^{-4}$ \\ $10^{-1}$ \end{tabular}&\begin{tabular}{@{}c@{}} $10$ \\ $8$ \\ $5$ \end{tabular}&\begin{tabular}{@{}c@{}} $48$ \\ $33$ \\ $17$ \end{tabular} & 0.99 & \begin{tabular}{@{}c@{}} $10^{-7}$ \\ $10^{-4}$ \\ $10^{-1}$ \end{tabular}&\begin{tabular}{@{}c@{}} $6$ \\ $5$ \\ $3$ \end{tabular}&\begin{tabular}{@{}c@{}} $44$ \\ $30$ \\ $16$ \end{tabular}\\

\hline
\end{tabular}
    \caption{Comparison of the iteration number $N$ to achieve  condition~\eqref{f9} and the estimate $N_*$ from  Theorem~\ref{thm-main}.}
    \label{tab:k_N}
\end{table}

In Table~\ref{tab:k_N} we can see that in all cases $N<N_*$. It means that  stopping condition~\eqref{f9} is reached earlier than the theoretical estimate of the number of iterations $N_*$ justified using PL-condition (see Theorem~\ref{thm-main}).  Also, we can note that the method converges much faster than the stated estimate for large values of $\mu$. At the same time, for small values of $\mu$, the value of $N_*$ exceeds $N$ by at most 2.5 times. For the other types of the noise of the gradient, a similar picture is observed.

Now, we compare the results of the Gradient Descent with a constant step-size~\eqref{f5} and the proposed Gradient Descent method with an adaptive step-size (Algorithm~\ref{adapt_gd}) when using stopping criterion~\eqref{f9}. In Tables~\ref{tab:qp_NT} and~\ref{tab:qp_dist}, there are  presented the results of the experiments for the quadratic function~\eqref{Q_function}. The experiments were carried out for  the uniformly distributed noise $v(x)$ on the sphere. In these experiments, the inexactness $\delta=16\Delta^2$  in the function was taken. Note that in this case, the correlation of inexactness satisfies the condition of Theorem~\ref{theorem:adaptL_inexactf}.

\begin{table}[ht]
    \centering
    \begin{tabular}{|c|c|c|c|c|c|c|c|}
      \hline
            &           &\multicolumn{2}{|c|}{Constant} & \multicolumn{2}{|c|}{Adaptive $L$} \\
            \hline
      $\mu$ & $\Delta$ & Iters & Time, ms & Iters & Time, ms \\
      \hline
0.01 & \begin{tabular}{@{}c@{}} $10^{-7}$ \\ $10^{-4}$ \\ $10^{-1}$ \end{tabular}&\begin{tabular}{@{}c@{}} $1525$ \\ $837$ \\ $158$ \end{tabular}&\begin{tabular}{@{}c@{}} $139.02$ \\ $76.72$ \\ $14.88$ \end{tabular}&\begin{tabular}{@{}c@{}} $668$ \\ $421$ \\ $75$ \end{tabular}&\begin{tabular}{@{}c@{}} $243.45$ \\ $151.83$ \\ $24.67$ \end{tabular}\\
\hline
0.1 & \begin{tabular}{@{}c@{}} $10^{-7}$ \\ $10^{-4}$ \\ $10^{-1}$ \end{tabular}&\begin{tabular}{@{}c@{}} $169$ \\ $104$ \\ $42$ \end{tabular}&\begin{tabular}{@{}c@{}} $15.84$ \\ $10.00$ \\ $4.40$ \end{tabular}&\begin{tabular}{@{}c@{}} $72$ \\ $46$ \\ $19$ \end{tabular}&\begin{tabular}{@{}c@{}} $26.40$ \\ $16.39$ \\ $6.62$ \end{tabular}\\
\hline
0.99 & \begin{tabular}{@{}c@{}} $10^{-7}$ \\ $10^{-4}$ \\ $10^{-1}$ \end{tabular}&\begin{tabular}{@{}c@{}} $6$ \\ $5$ \\ $3$ \end{tabular}&\begin{tabular}{@{}c@{}} $1.01$ \\ $0.92$ \\ $0.55$ \end{tabular}&\begin{tabular}{@{}c@{}} $6$ \\ $5$ \\ $3$ \end{tabular}&\begin{tabular}{@{}c@{}} $2.16$ \\ $1.75$ \\ $1.12$ \end{tabular}\\
\hline
\end{tabular}
    \caption{Comparison of the running time of the algorithms and number of iterations to achieve the accuracy $\|\widetilde{\nabla} f(x)\|\leqslant \sqrt{6}\Delta$ for the quadratic problem.}
    \label{tab:qp_NT}
\end{table}

From Table~\ref{tab:qp_NT}, we can see that the adaptive method is inferior in real time to the Gradient Descent for all parameters $\mu$ and $\Delta$. However, it needs a smaller number of iterations for big values of $\mu$.

\subsection{The problem of minimizing the logistic regression function}

Now let us  check the work of the proposed stopping criterion in the case when it is rather difficult to estimate the constant $\mu$ of the function which satisfies PL-condition. In this case, we will not be able to use estimate~\eqref{equat_estimN*}. This situation has been discussed in Remark~\ref{rem-notmu}. The detailed experiments presented in~\ref{sec:app_exp_logreg}.

However, we note that, as shown by the previous experiment, condition~\eqref{f9} can be achieved in a significantly smaller number of steps compared to the theoretical estimate of the number of iterations $N_*$ from Theorem~\ref{thm-main}.
We will consider the following optimization problem associated with logistic regression
    $f(x) = \frac{1}{m}\sum\limits_{i=1}^m\log \left(1 + \exp\left(-y_i \langle w_i, x\rangle\right)\right ),$
where $y=\left(y_1,\dots, y_m\right)^\top \in [-1,1]^m$ is the feasible variable vector, $W=[w_1\dots w_m]\in \mathbb{R}^ {n\times m}$ is the feature matrix, where the vector $w_i\in\mathbb{R}^n$ is from the same space as the optimized weight vector $w$.

Note that this problem may not have a finite solution in the general case. So we will create such an artificial data set that there is a finite vector $x^*$ minimizing the given function. The details of data generation is presented in Appendix~\ref{sec:app_exp_logreg}.

\begin{figure}[ht]
\centering
\includegraphics[width=0.5\linewidth]{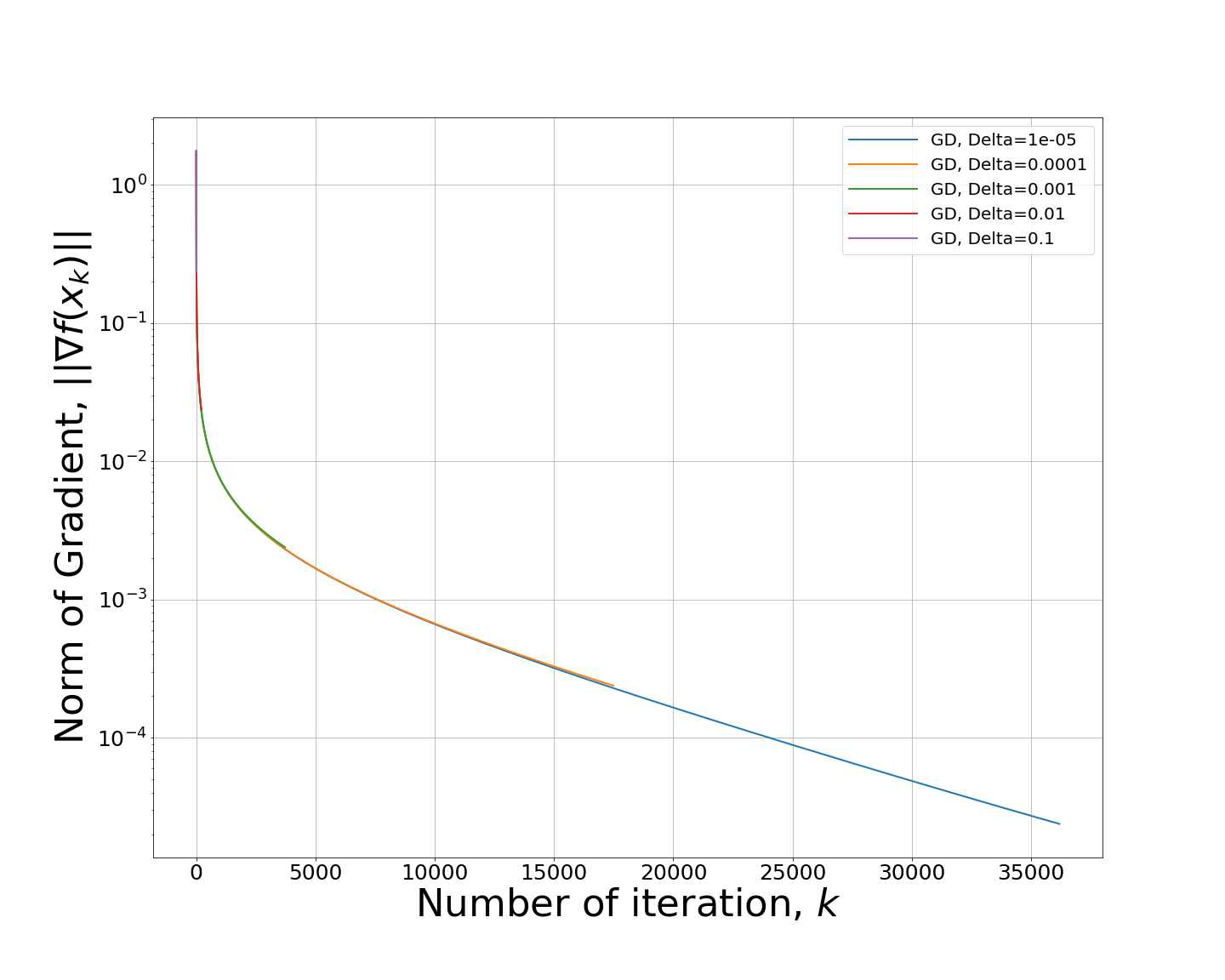}
\caption{The rate of convergence of the gradient method in the gradient norm for different values of $\Delta$ for the problem of minimizing logistic regression using stopping criterion~\eqref{f9} for the constant $v$.} \label{fig:grad4}
\end{figure}

In the conducted experiments, we chose $n=200, m=700$ and $k=10< \min\left(n,\frac{m}{2}\right)$.
We consider in this section the case of constant inexactness. From Fig.~\ref{fig:grad4} it can be seen that the trajectories of the method are not the same. Moreover, adding inexactness slows down the convergence. On the other hand, the trajectories have become more similar compared to the case of the inexactness directed along the minus of the gradient (see Appendix~\ref{sec:app_exp_logreg}).

\begin{figure}[ht]
\vspace{-4ex} \centering \subfigure[]{
\includegraphics[width=0.45\linewidth]{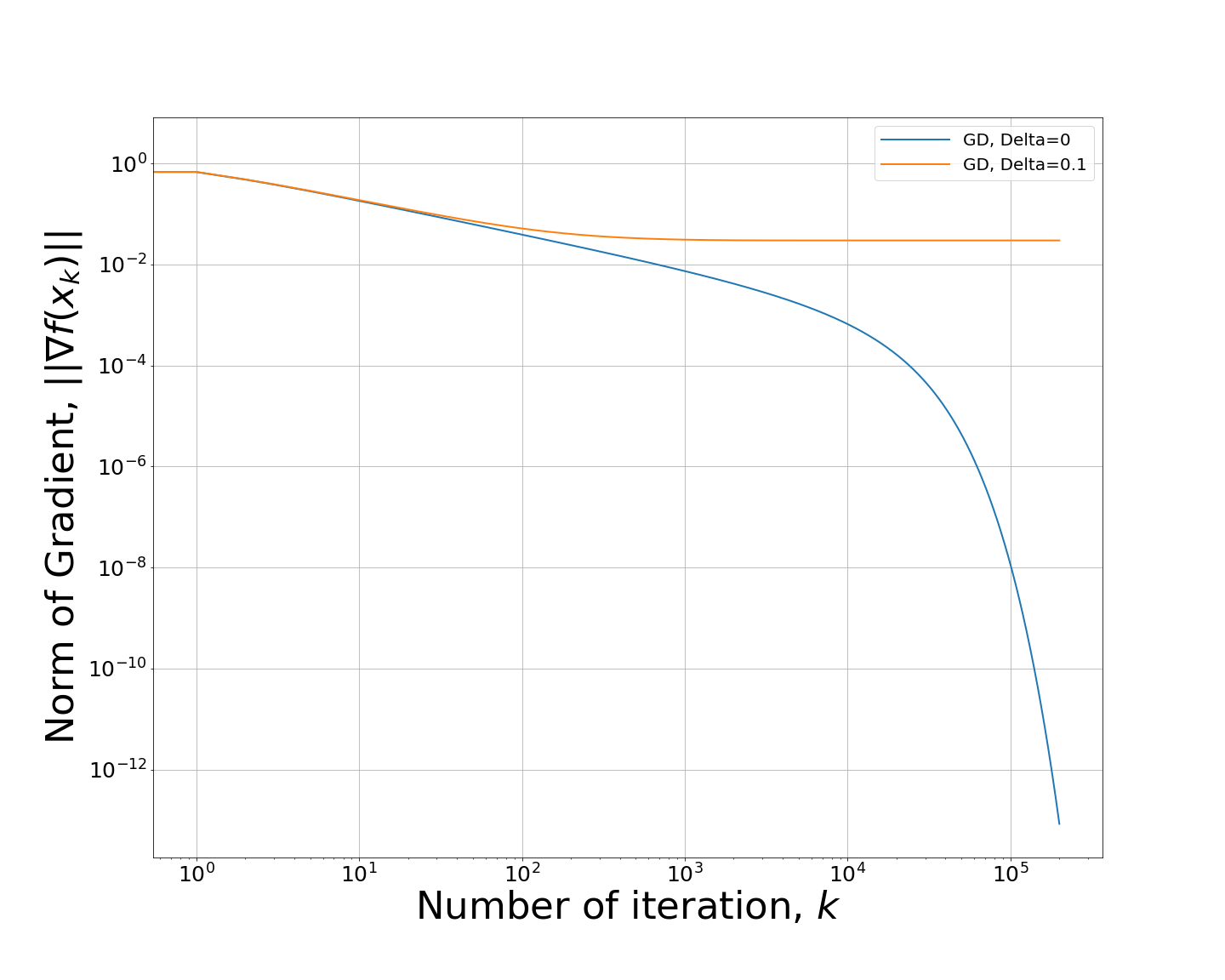}  \label{fig:grad4_delta} }
\hspace{4ex}
\subfigure[]{
\includegraphics[width=0.45\linewidth]{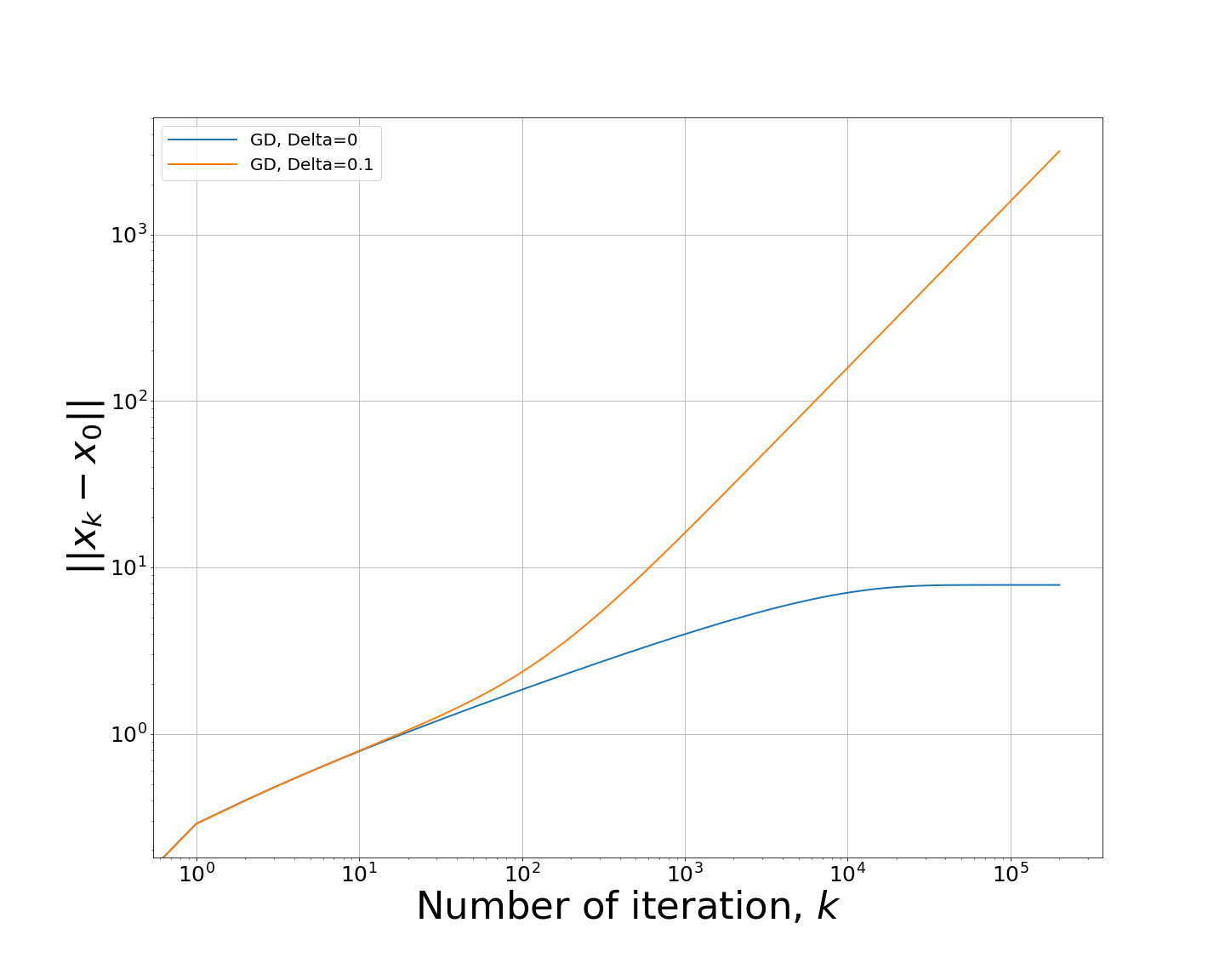}  \label{fig:norm_x4_delta} }
\caption{The results of the gradient method with respect to the norm of the gradient without using the stopping criterion for $\Delta=0.1$ for the problem of logistic regression minimization for the constant inaccuracy $\Delta v$. \subref{fig:grad4_delta} The convergence rate with respect to the norm of the gradient; \subref{fig:norm_x4_delta} the distance from the starting point to $x_k$.} \label{fig:dtype4}
\end{figure}

However, in this case, in Fig.~\ref{fig:norm_x4_delta} it can be seen that without using the stopping criterion, the distance $\|x_k - x_0\|$ grows rather quickly. Thus, in the case of randomly generated gradient noise after $10^5$ iterations of Gradient Descent method~\eqref{f5} the distance was $1.25$ times larger compared to the result without noise in the gradient. At the same time, in the case of a constant gradient specification error, these values differ by more than two orders of magnitude  (see Appendix~\ref{sec:app_exp_logreg}).

\subsection{Some experiments with the Rosenbrock-type function}
In this subsection, we describe results of our investigation of the behavior of the proposed adaptive Algorithm~\ref{adapt_gd} for  some non-convex problems. The detailes are presented in Appendix~\ref{sec:app_exp_rosenbrock} and~\ref{sec:app_exp_ns}. Firstly, we considered the well-known two-dimensional Rosenbrock function
$f(x_1, x_2) = 100(x_2-(x_1)^2)^2 + (x_1-1)^2.$

This function is not convex, and it satisfies the Lipschitz condition for the gradient only locally. Indeed, if we consider the line $x_2=0$, then we get  $f(x_1,0)=100x_1^4+(x_1-1)^2$. The gradient of this function does not satisfy the Lipschitz condition. On the other hand, the Rosenbrock function satisfies locally PL-condition.

In the conducted experiments, we will vary the value of the parameter $\Delta$ and take $\delta=\Delta^2.$ In Table~\ref{tab:rosenbrock2} in Appendix~\ref{sec:app_exp_rosenbrock}, we show the results for different types of noise. As previously, from the results presented  in Table~\ref{tab:rosenbrock2}, we can see that the number of required iterations increases with decreasing $\Delta$ (which also tightens the stopping condition). Moreover, it increases logarithmically, which coincides with the results of Theorem~\ref{theorem:adaptL_inexactf}. We can also note that the resulting distance from the starting point $x_0$ to the last point does not exceed the distance from the starting  point $x_0$  to the nearest optimal one $x_* = (1,1)$ everywhere. In addition, for all considered types of the gradient error (noise), a comparable convergence rate is observed according to the number of iterations until stopping criterion~\eqref{stop_cond_adaptGD} is satisfied, and to the running time for the corresponding values of $\Delta$.

Further, let us consider a system of nonlinear equations $g(x)=0,$ where $g_1=\frac{1}{2}(x_1-1), g_i=x_{i}-2x_{i-1}^2+1, i =\overline{2,n}.$ The problem of solving this system is equivalent to minimizing the following Nesterov-Skokov function (see~\cite{NestSkok})
\begin{equation}
    \label{NSfunction}
    f(x)=\frac{1}{4}\left(1-x_1\right)^2+\sum\limits_{i=1}^{n-1} \left(x_{i+1}-2x_i^2+1\right)^2.
\end{equation}
This function is analogous to the Rosenbrock function. It is also non-convex and satisfies the Lipschitz gradient condition only locally. Also,  function~\eqref{NSfunction} has a global minimum at the point $(1,1\dots 1,1)^\top$ and an optimal value $f^*=0$. Moreover, this function locally satisfies PL-condition (see the proof in Appendix~\ref{app:NS_function}).

As it was seen from the results of the previous experiments, our proposed stopping criterion~\eqref{stop_cond_adaptGD} of Algorithm~\ref{adapt_gd} can work equally well for all considered types of noise in the gradient. In the current experiments for the Nesterov-Skokov function, we used the random noise of the gradient which is uniformly distributed  on the sphere. For the experiments, the starting point is $(-1,1,\dots 1, 1)^\top$ and therefore $\|x_0-x_*\|=2$. We will vary the value of the inexactness $\Delta$ and the dimension of the problem $n$.

Table~\ref{tab:NestSkok} in Appendix~\ref{sec:app_exp_ns} shows the results of the adaptive gradient method~\ref{adapt_gd} for the Nesterov-Skokov function~\eqref{NSfunction}. Firstly, we see that as the dimension of $n$ increases, the difference between the required time to solve the problem for different $\Delta$ grows significantly. Secondly, for different $n$ with the same $\Delta$, the method converges to a solution with significantly different accuracy. We can also note that $\|x_N-x_0\|$ exceeds $\|x_0-x_*\|$ by at most 2 times. Moreover, significant upward deviations are observed for the cases when numerous iterations are made ($n=5$ and $\Delta=10^{-4}, 10^{-3}$). It can also be noted that even for sufficiently small values of the norm of the gradient, the accuracy by the function turns out to be quite low (which is typical for the Nesterov-Skokov function).

\section{Conclusion}

This paper studies stopping criterions for the gradient method with an inexact gradient. The authors focus on the case of non-convex functions. The paper presents a stopping criterion that finds a compromise between the accuracy of the obtained point and the distance to the starting point. Moreover, it is shown that the method moves away from the starting point to a distance comparable to the distance to the nearest solution if the function satisfies PL-condition.

Besides, the paper considers the cases of a constant and adaptive step size in the gradient methods. For both cases, we present theoretical analysis and the number of iterations required to approach the stopping criterion or to find the point with the required quality.

In addition, the paper contains numerical experiments demonstrating the work of the stopping criterion. In particular, there are experiments on a quadratic function (convex, but not strongly convex), demonstrating the stopping criterion to be approached faster than the theoretical estimation of the iteration number $N_*$. Also, we present experiments on the problem of logistic regression where the objective function is convex and meets PL-condition only locally. The proposed stopping criterion on this function stops the growth of the distance $\|x_k-x_0\|$. Moreover, we present experiments on non-convex functions: the Rosenbrock function and its multidimensional generalization, which is the Nesterov-Skokov function. The first function demonstrates that our stopping criterion works with general types of inexactness. The second function demonstrates that even a small inexactness can lead to quite a high value of function and this value cannot be improved. Also, we demonstrate that for some noises, the gradient method can move away quite far on the Nesterov-Skokov function without a stopping criterion.

\phantomsection\addcontentsline{toc}{section}{Acknowledgements}
\section*{Acknowledgements}

The research by F.\,Stonyakin in Section 2.1 was supported  by the strategic academic leadership program "Priority 2030" (Agreement  075-02-2021-1316, 30.09.2021). The research by B.\,Polyak  in Section 1 was supported by the Russian Science Foundation (project No. 21-71-30005).

\newpage

\phantomsection\addcontentsline{toc}{section}{\refname}
\bibliographystyle{spmpsci}
\bibliography{biblio_new}

\newpage

\appendix

\section{The proof of Theorem \ref{theorem:adaptL_inexactf}}
\label{app:adaptL_inexactf_proof}
In the case, if at some $k$-th iteration of Algorithm \ref{adapt_gd} the condition \eqref{stop_cond_adaptGD} is satisfied, then according to the PL-condition we obtain that
$$
    f(x_k)-f^*\leq \frac{5\Delta^2}{\mu}.
$$

Let us study the estimation of the quality of the output point of Algorithm \ref{adapt_gd} under conditions when \eqref{stop_cond_adaptGD} is not satisfied. Note that for each iteration $k\geq 1$  condition \eqref{cond_it} is satisfied, i.e. 
$$
    \widetilde{f}(x_{k+1})\leq \widetilde{f}(x_k)+\langle\widetilde{\nabla} f(x_k), x_{k+1}-x_k\rangle + L_k\|x_{k+1}-x_k\|^2 + \frac{\Delta^2}{2L_k} +2\delta.
$$
Then by using  condition \eqref{inexactf_cond}, we get
$$
    f(x_{k+1})\leq f(x_k)+\langle\widetilde{\nabla} f(x_k), x_{k+1}-x_k\rangle + L_k\|x_{k+1}-x_k\|^2 + \frac{\Delta^2}{2L_k} +4\delta.
$$

Moreover, taking into account \eqref{x_k}, we have the following inequality
\begin{equation}\label{k_inequality}
    f(x_{k+1})-f(x_k)\leq \frac{\Delta^2}{2L_k} - \frac{1}{4L_k}\|\widetilde{\nabla}f(x_k)\|^2 +2\delta.
\end{equation}

This condition, together with the relation for the inexactness $8L_k\delta\leq\Delta^2$, tells us that if $\|\widetilde{\nabla}f(x_k)\|\geq C\Delta$ for $C>\sqrt{3}$, then $f(x_{k+1}) < f(x_k)$ is guaranteed and the method converges to the minimum. For definiteness, we take $C=2$. 
 Then if $\|\widetilde{\nabla} f(x_{k})\|>2\Delta $ and \eqref{k_inequality} holds for every $k=0,1,\ldots,N-1$, then
$$
    f(x_{0})-f(x_{N})=\sum_{k=0}^{N-1}(f(x_{k})-f(x_{k+1}))>\frac{\Delta^{2}}{2}\sum_{k=0}^{N-1}\frac{1}{L_k}-2\delta N\geq \frac{N\Delta^{2}}{4L}-4\delta N,
$$
i.e. at $16L\delta < \Delta^2$, we have
\begin{equation}\label{N_delta}
    N < \frac{2L}{\Delta^{2}-16L\delta}(f(x_{0})-f^{*}),
\end{equation}
which indicates the end of the process. 

On the other hand, from estimate \eqref{k_inequality} one can get an estimate for the function residual at the $k$-th iteration. Using the PL-condition, we get
$$ 
    f(x_{k+1})-f(x_k)\leq \frac{\Delta^2}{2L_k} - \frac{\mu}{4L_k}(f(x_k)-f^*) +4\delta,
$$
whence
\begin{align*}
    f(x_{k+1})-f^*&\leq \left(1-\frac{\mu}{4L_k}\right)(f(x_k)-f^*) + \frac{\Delta^2}{2L_k}+2\delta \\
    &\leq \prod_{j=0}^k \left(1-\frac{\mu}{4L_j}\right)(f(x_0)-f^*) + \frac{\Delta^2}{2}\left(\frac{1}{L_k}+\sum\limits_{j=0}^{k-1} \frac{1}{L_{j}} \prod_{i=j+1}^{k} \left(1-\frac{\mu}{4L_i}\right)\right) +\\
    & \qquad  + 4\delta\left(\sum\limits_{j=0}^{k-1}\prod_{i=j+1}^{k} \left(1-\frac{\mu}{4L_i}\right)\right).
\end{align*}


Let us estimate the second term. We denote by  $S_k = \left(\frac{1}{L_k}+\sum\limits_{j=0}^{k-1} \frac{1}{L_{j}} \prod_{i=1}^{k-j} \left(1-\frac{\mu}{4L_i}\right)\right)$. 
Note that $S_k$ with $k\geq 1$ satisfies the recursive formula
$S_k=\frac{1}{L_k}+S_{k-1}\left(1-\frac{\mu}{4L_k}\right).$  Let us consider two cases. In the first case, if $\mu S_{k-1}\geq 4$, then $S_k=\frac{4-\mu S_{k-1}}{4L_k}+S_{k-1}\leq S_{k-1}$. In the second case, for $\mu S_{k-1}< 4$ we get that $S_{k-1}<\frac{4}{\mu}$ and $S_k< \frac{1}{L_k}+\frac{4}{\mu}\left(1-\frac{\mu}{4L_k}\right)=\frac{4}{\mu}$. Thus $S_k\leq \max\left(S_{k-1}, \frac{4}{\mu}\right)$. By sequentially expanding, we obtain the estimate 
$$
    S_k\leq \max\left(S_0, \frac{4}{\mu}\right)=\max\left(\frac{1}{L_0}, \frac{4}{\mu}\right)\leq \frac{4}{\mu}.
$$

We also estimate the third term using $L_j\leq 2L$ as
$$
    \sum\limits_{j=0}^{k-1}\prod_{i=j+1}^{k} \left(1-\frac{\mu}{4L_i}\right)\leq \sum\limits_{j=0}^{k-1}\left(1-\frac{\mu}{8L}\right)^j \leq \frac{8L}{\mu}.
$$

Thus, we obtain the estimate
\begin{equation}
\label{exact_fk}
    f(x_{k+1})-f^*\leq  \prod_{j=0}^k \left(1-\frac{\mu}{4L_j}\right)(f(x_0)-f^*) + \frac{2\Delta^2}{\mu} + \frac{32L}{\mu}\delta.
\end{equation}

The estimate of \eqref{exact_fk} depends on how the algorithm works. Using the inequality $L_k\leq 2L$, we obtain a result about the convergence which is determined only by the parameters of the function
\begin{equation}
\label{inexact_fk1}
    f(x_{k+1})-f^*\leq  \left(1-\frac{\mu}{8L}\right)^{k+1}(f(x_0)-f^*)  + \frac{2\Delta^2}{\mu}  + \frac{32L}{\mu}\delta.
\end{equation}
Next, we use the relation for the inexactness $16L\delta\leq\Delta^2$ and  estimate \eqref{inexact_fk1}:
\begin{equation}
\label{inexact_fk}
    f(x_{k+1})-f^*\leq  \left(1-\frac{\mu}{8L}\right)^{k+1}(f(x_0)-f^*)  + \frac{4\Delta^2}{\mu}.
\end{equation}

Let us estimate the distance $\|x_0-x_N\|$ in the same way as it was done for the gradient descent method with a constant step. Using  inequalities \eqref{k_inequality}, \eqref{inexact_fk} and and taking into account that  $\widetilde{\nabla}f(x_k)=2L_k(x_{k+1}-x_k)$, we get the following estimate
\begin{align*}
    \|x_{k+1}-x_k\|^2&=\frac{\|\widetilde{\nabla}f(x_k)\|^2}{4L_k^2} \leq \frac{\Delta^2}{2L_k^2} + \frac{f(x_{k+1})-f(x_k)}{L_k}\\
    &\leq \frac{\Delta^2}{2L_k^2} + \frac{f(x_{k+1})-f^*}{L_k}\\
    &\leq \frac{\Delta^2}{2L_k^2}+ \frac{4\Delta^2}{\mu L_k} + \frac{1}{L_k}\left(1-\frac{\mu}{8L}\right)^{k}(f(x_0)-f(x^*))\\
    &\leq \frac{4\Delta^2}{2L_{\min}^2} + \frac{\Delta^2}{\mu L_{\min}}+ \frac{1}{L_{\min}}\left(1-\frac{\mu}{8L}\right)^{k}(f(x_0)-f(x^*)).
\end{align*}

After summing these inequalities over $k$ from $0$ to $N-1$, we get the following result
$$
    \|x_{k+1}-x_k\| \leq \Delta\sqrt{\frac{1}{2L_{\min}^2} + \frac{4}{\mu L_{\min}}} +  \left(1-\frac{\mu}{8L}\right)^{\frac{k}{2}} \sqrt{\frac{f(x_0)-f(x^*)}{L_{\min}}}.
$$
Thus,
\begin{equation} \label{distN_adaptL}
    \|x_N-x_0\|\leq \sum\limits_{k=0}^{N-1}\|x_{k+1}-x_k\|\leq N\Delta \sqrt{\frac{1}{2L_{\min}^2} + \frac{4}{\mu L_{\min}}} + 16\sqrt{\frac{L}{L_{\min}}}\frac{\sqrt{L(f(x_0)-f^*)}}{\mu}.
\end{equation}
We can estimate the number of iterations of Algorithm \ref{adapt_gd} from  the inequality \eqref{inexact_fk} as follows
$$
    N\leq \left\lceil\frac{8L}{\mu}\log \frac{\mu (f(x_0)-f^*)}{2\Delta^2}\right\rceil.
$$

Let us introduce the notation $\gamma = \frac{L}{L_{\min}}$. Then we can estimate the final estimation of the distance from the starting point $x_0$ to the current one $x_N$ as follows
$$
    \|x_N-x_0\|\leq 8\frac{\Delta}{\mu} \sqrt{\frac{1}{2}\gamma^2  + 4\gamma\frac{L}{\mu }} \log \frac{\mu (f(x_0)-f^*)}{2\Delta^2} + 16\frac{\sqrt{\gamma L(f(x_0)-f^*)}}{\mu}.
$$

Note that the factor $\gamma=\frac{L}{L_{\min}}$ appeared in this estimate, which depends on an unknown constant and the parameter of the algorithm $L_{\min}$. In the worst case, we have $\gamma = \frac{4L}{\mu}$.

Let us estimate the number of additional calculations of the value of the function and operations of the form \eqref{x_k} of the adaptive gradient method \ref{adapt_gd} in comparison with the gradient descent with a constant step-size \eqref{f5}. Let $i_k$ computations of the inexact gradient be made at the $k$-th step. Then $2^{i_1-1}=\frac{2L_k}{L_{k-1}}$. Then we note that
$$
    \prod_{k=1}^N 2^{i_k-1} = \prod_{k=1}^N \frac{2L_k}{L_{k-1}}= 2^N \frac{L_N}{L_0}.
$$

As mentioned above, it is true that $L_N \leq 2L$. Then the total number of additional function evaluations and steps of the form \eqref{x_k} $I(N)=\sum\limits_{i=1}^N i_k$ is estimated from above as follows
$$
    2^{I(N)-2N} \leq \frac{2L}{L_0}.
$$
Therefore,
$$
    I(N)\leq 2N + \log \frac{2L}{L_0}.
$$

\section{Numerical Experiments}

\subsection{The quadratic form minimization problem}
\label{sec:app_exp_qp}

In this section, we present additional experiments for quadratic problem presented in Section \ref{sec:exp_qp}.

\begin{figure}[ht!]
\centering
\includegraphics[width=0.4\linewidth]{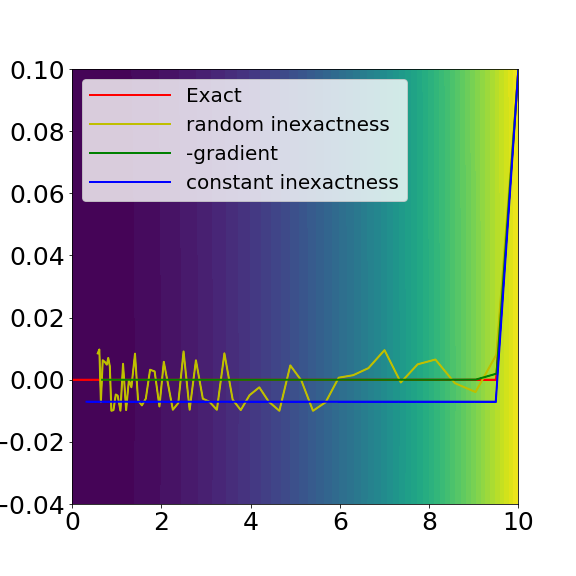}
\caption{Convergence of the gradient method with stopping criterion \eqref{f9} in the two-dimensional case for different $v$: zero (exact), randomly generated at each iteration (random inexactness), co-directional with minus gradient (Antigradient), and constant (constant inexactness).} \label{fig:lines}
\end{figure}

Now, let us consider a quadratic form in dimension $n=2$ with coefficients 0.05 and 1. In this case, we can draw the convergence trajectory of the method for different types of noise (see Fig. \ref{fig:lines}). We chose the starting point  $x_0=(10,0.1)^\top$ and noise level  $\Delta=10^{-2}.$ In the case of constant inexactness, the vector $v=(1,1)^\top$ was taken.

We can see that, if the inexact gradient is collinear to the exact one, then the trajectories coincide almost completely for the gradient method with exact and inexact gradients. At the same time, we do not see changes in the coordinate $x_k$, on which the function almost does not depend. For the constant noise, we observe a constant displacement of the trajectory from the path of the inexact gradient descent method.

Now we compare the results of the gradient descent with a constant step-size \eqref{f5} and the proposed gradient descent method with an adaptive step-size (Algorithm \ref{adapt_gd}) when we use stopping criterion \eqref{f9}. In Tables \ref{tab:qp_NT} and \ref{tab:qp_dist}, there were  presented the results of the experiments for the quadratic function \eqref{Q_function}. The experiments were carried out for  the uniformly distributed noise $v(x)$ on the sphere. In these experiments, the inexactness $\delta=16\Delta^2$  in the function was taken. Note that in this case, the correlation of inexactness satisfies the condition of Theorem \ref{theorem:adaptL_inexactf}.

From Table \ref{tab:qp_NT}, we can see that the adaptive method is inferior in real time to the gradient descent for all parameters $\mu$ and $\Delta$. However, it needs a smaller number of iterations for small $\mu$.

\begin{table}[ht]
    \centering
    \begin{tabular}{|c|c|c|c|c|c|c|c|}
      \hline
            &   &\multicolumn{3}{|c|}{Constant} & \multicolumn{3}{|c|}{Adaptive $L$}\\
            \hline
      $\mu$ &  $\Delta$ & $\|x_N-x_0\|$ & $\frac{\|\nabla f(x_N)\|}{\Delta}$ & $f(x_N)-f^*$& $\|x_N-x_0\|$ & $\frac{\|\nabla f(x_N)\|}{\Delta}$ & $f(x_N)-f^*$\\
      \hline
0.01 & \begin{tabular}{@{}c@{}} $10^{-7}$ \\ $10^{-4}$ \\ $10^{-1}$ \end{tabular}&\begin{tabular}{@{}c@{}} $948.7$ \\ $948.7$ \\ $946.2$ \end{tabular}&\begin{tabular}{@{}c@{}} $2.34$ \\ $2.37$ \\ $2.39$ \end{tabular}&\begin{tabular}{@{}c@{}} $0.25 \cdot 10^{-11}$ \\ $0.26 \cdot 10^{-5}$ \\ $2.50$ \end{tabular}&\begin{tabular}{@{}c@{}} $948.7$ \\ $948.7$ \\ $946.5$ \end{tabular}&\begin{tabular}{@{}c@{}} $2.01$ \\ $2.04$ \\ $2.30$ \end{tabular}&\begin{tabular}{@{}c@{}} $0.18 \cdot 10^{-11}$ \\ $0.19 \cdot 10^{-5}$ \\ $1.93$ \end{tabular}\\
\hline
0.1 & \begin{tabular}{@{}c@{}} $10^{-7}$ \\ $10^{-4}$ \\ $10^{-1}$ \end{tabular}&\begin{tabular}{@{}c@{}} $948.7$ \\ $948.7$ \\ $948.3$ \end{tabular}&\begin{tabular}{@{}c@{}} $2.17$ \\ $1.95$ \\ $2.15$ \end{tabular}&\begin{tabular}{@{}c@{}} $0.22 \cdot 10^{-12}$ \\ $0.17 \cdot 10^{-6}$ \\ $0.20$ \end{tabular}&\begin{tabular}{@{}c@{}} $948.7$ \\ $948.7$ \\ $948.4$ \end{tabular}&\begin{tabular}{@{}c@{}} $1.90$ \\ $2.14$ \\ $1.80$ \end{tabular}&\begin{tabular}{@{}c@{}} $0.14 \cdot 10^{-12}$ \\ $0.19 \cdot 10^{-6}$ \\ $0.14$ \end{tabular}\\
\hline
0.9 & \begin{tabular}{@{}c@{}} $10^{-7}$ \\ $10^{-4}$ \\ $10^{-1}$ \end{tabular}&\begin{tabular}{@{}c@{}} $948.7$ \\ $948.7$ \\ $948.7$ \end{tabular}&\begin{tabular}{@{}c@{}} $0.93$ \\ $0.86$ \\ $0.99$ \end{tabular}&\begin{tabular}{@{}c@{}} $0.46 \cdot 10^{-14}$ \\ $0.39 \cdot 10^{-8}$ \\ $0.52 \cdot 10^{-2}$ \end{tabular}&\begin{tabular}{@{}c@{}} $948.7$ \\ $948.7$ \\ $948.7$ \end{tabular}&\begin{tabular}{@{}c@{}} $0.92$ \\ $0.95$ \\ $0.88$ \end{tabular}&\begin{tabular}{@{}c@{}} $0.45 \cdot 10^{-14}$ \\ $0.48 \cdot 10^{-8}$ \\ $0.41 \cdot 10^{-2}$ \end{tabular}\\
\hline
0.99 & \begin{tabular}{@{}c@{}} $10^{-7}$ \\ $10^{-4}$ \\ $10^{-1}$ \end{tabular}&\begin{tabular}{@{}c@{}} $948.7$ \\ $948.7$ \\ $948.7$ \end{tabular}&\begin{tabular}{@{}c@{}} $0.99$ \\ $0.94$ \\ $1.02$ \end{tabular}&\begin{tabular}{@{}c@{}} $0.49 \cdot 10^{-14}$ \\ $0.44 \cdot 10^{-8}$ \\ $0.52 \cdot 10^{-2}$ \end{tabular}&\begin{tabular}{@{}c@{}} $948.7$ \\ $948.7$ \\ $948.6$ \end{tabular}&\begin{tabular}{@{}c@{}} $0.94$ \\ $0.93$ \\ $0.99$ \end{tabular}&\begin{tabular}{@{}c@{}} $0.44 \cdot 10^{-14}$ \\ $0.43 \cdot 10^{-8}$ \\ $0.49 \cdot 10^{-2}$ \end{tabular}\\
\hline
\end{tabular}
    \caption{Comparison of algorithms in terms of the achieved accuracy in terms of the gradient norm and the distance from the start point to the last point. The distance from the starting point to the nearest optimal one is $948.7$.}
    \label{tab:qp_dist}
\end{table}

Note that according to the results presented in Table \ref{tab:qp_dist}, the compared methods lead to the achievement of approximately the same quality of the approximate solution. In this case, the trajectories of the methods are moved away approximately equally from the starting point. In this case, the distance $\|x_N - x_0\|$ is approximately equal to the distance $x_0$ to the nearest exact solution $x_*$. As we can see, the achieved accuracy is no less than $\frac{7\Delta^2}{\mu}$.

\subsection{The problem of minimizing the logistic regression function}
\label{sec:app_exp_logreg}

Now let us  check the work of the proposed stopping criterion in the case when it is rather difficult to estimate the constant $\mu$ of the function which satisfies the PL-condition. In this case, we will not be able to use estimate \eqref{equat_estimN*}. This situation has been discussed in Remark \ref{rem-notmu}. 

However, we note that, as shown by the previous experiment, condition \eqref{f9} can be achieved in a significantly smaller number of steps compared to the theoretical estimate of the number of iterations $N_*$ from Theorem \ref{thm-main}.
We will consider the following optimization problem associated with logistic regression:
\begin{equation}
\label{log_regr}
    f(x) = \frac{1}{m}\sum\limits_{i=1}^m\log \left(1 + \exp\left(-y_i \langle w_i, x\rangle\right)\right ),
\end{equation}
where $y=\left(y_1,\dots, y_m\right)^\top \in [-1,1]^m$ is the feasible variable vector, $W=[w_1\dots w_m]\in \mathbb{R}^ {n\times m}$ is the feature matrix, where the vector $w_i\in\mathbb{R}^n$ is from the same space as the optimized weight vector $w$. 

Note that this problem may not have a finite solution in the general case. So we will create such an artificial data set that there is a finite vector $x^*$ minimizing the given function. To do this, we generate $W$ and $y$ as follows:
\begin{enumerate}
\item We construct  $k\leq \min\left(n,\frac{m}{2}\right)$ orthogonal vectors with the unit norm and combine them into a matrix $W_B\in\mathbb{R}^{n\times k}$.

\item Construct a matrix $\widetilde{W} = W_BV^\top \in\mathbb{R}^{n\times m-2k}, $ where $V\in\mathbb{R}^{m-2k\times k}$ is some random matrix that defines the expansion of the vectors from the matrix $W$ in the basis $W_B$.

\item Construct some vector $x_0$ and define new vectors $\widetilde{y}=\text{sign}(\widetilde{W}x_0)$, $y_1=\text{sign}(W_B x_0)$.

\item Define the feasible variable vector $y = [y_1|-y_1|\widetilde{y}] \in [-1,1]^m$.

\item Define the feature matrix $W = [W_B|W_B|\widetilde{W}]\in \mathbb{R}^{m\times n}$.
\end{enumerate}

\begin{proposition}
\label{prop:logreg}
For the function \eqref{log_regr}  with such data, the following statements hold:
\begin{enumerate}
    \item The function $f$ satisfies the PL-condition  on any compact set $K$;

    \item The function $f$ has a Lipschitz gradient with the constant $L=\frac{\lambda_{\max}(W^\top W)}{4m}$;

    \item There is a finite $x_*$,  where the objective function reaches its minimal value;

    \item The set of minimum points $X_*$ is unbounded if $k<n$.
\end{enumerate}
\end{proposition}

The proof of Proposition \ref{prop:logreg} is presented in Appendix \ref{proof_prop_logreg}.

In the conducted experiments, we chose $n=200, m=700$ and $k=10< \min\left(n,\frac{m}{2}\right)$. Accordingly, there are 700 objects with 200 features and a feature matrix of rank 10. Thus, the set of the  solutions for the minimization problem of the function $f$  with such data is non-empty and unbounded.

Let us apply the gradient descent method \eqref{f5} for various inexactness $\Delta$  in the gradient with stopping criterion \eqref{f9}. As an inexact gradient, we will use a gradient with random noise \textbf{Random}.

\begin{figure}[ht!]
\centering
\includegraphics[width=0.5\linewidth]{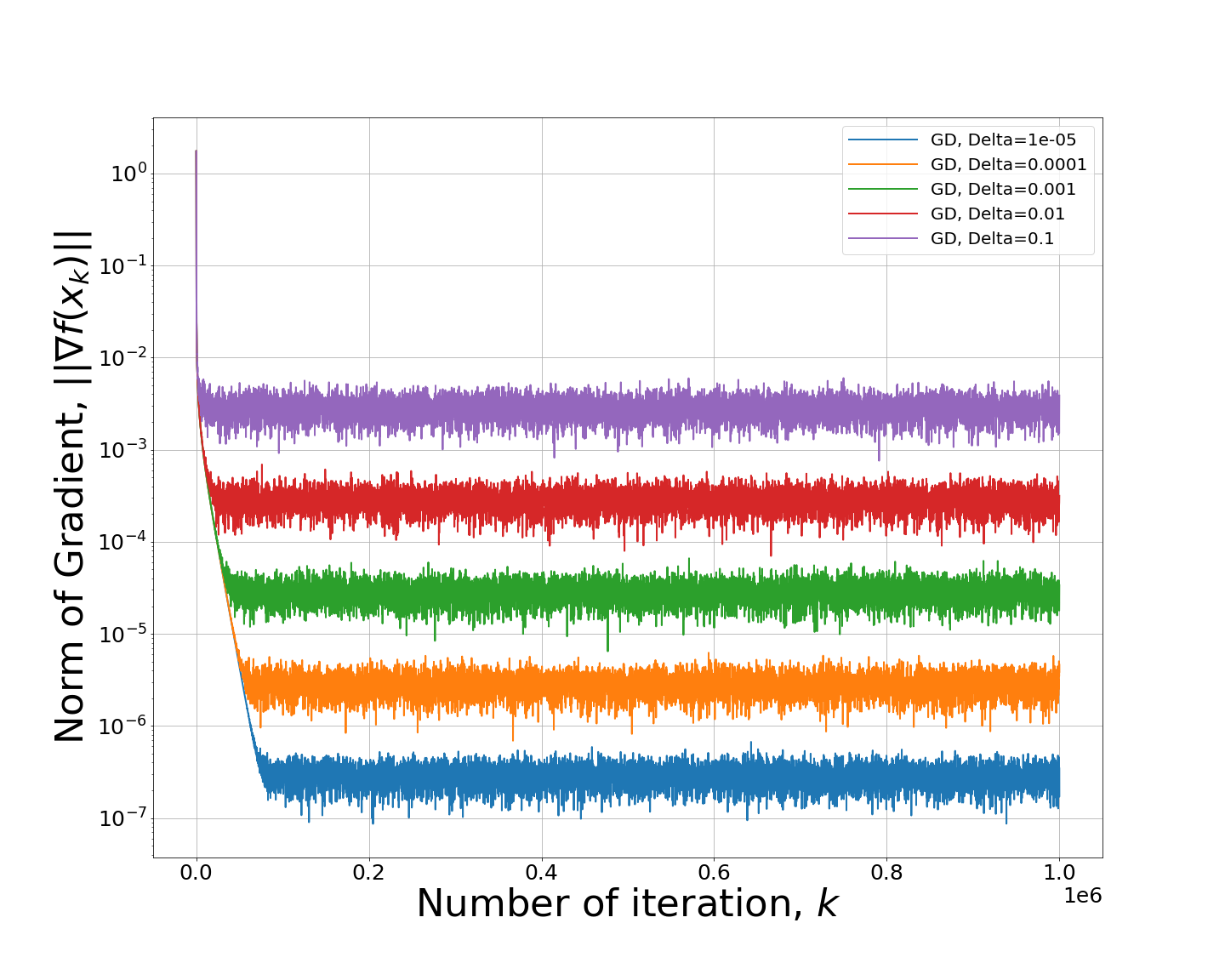}
\caption{The convergence rate of the gradient method with respect to the norm of the gradient for different values of the inexactness $\Delta$ for the problem of logistic regression minimization at the first $N=10^5$ iterations without using the stopping criterion.} \label{fig:grad_non_stop}
\end{figure}

In Fig. \ref{fig:grad_non_stop} there was shown the plot of the convergence of the gradient method with different levels of the noise  $\Delta$  without using the stopping criterion proposed  in this article. It can be seen that the method reaches points with a gradient norm of order $\Delta$, but it cannot converge closer.

\begin{figure}[ht]
\centering
\includegraphics[width=0.5\linewidth]{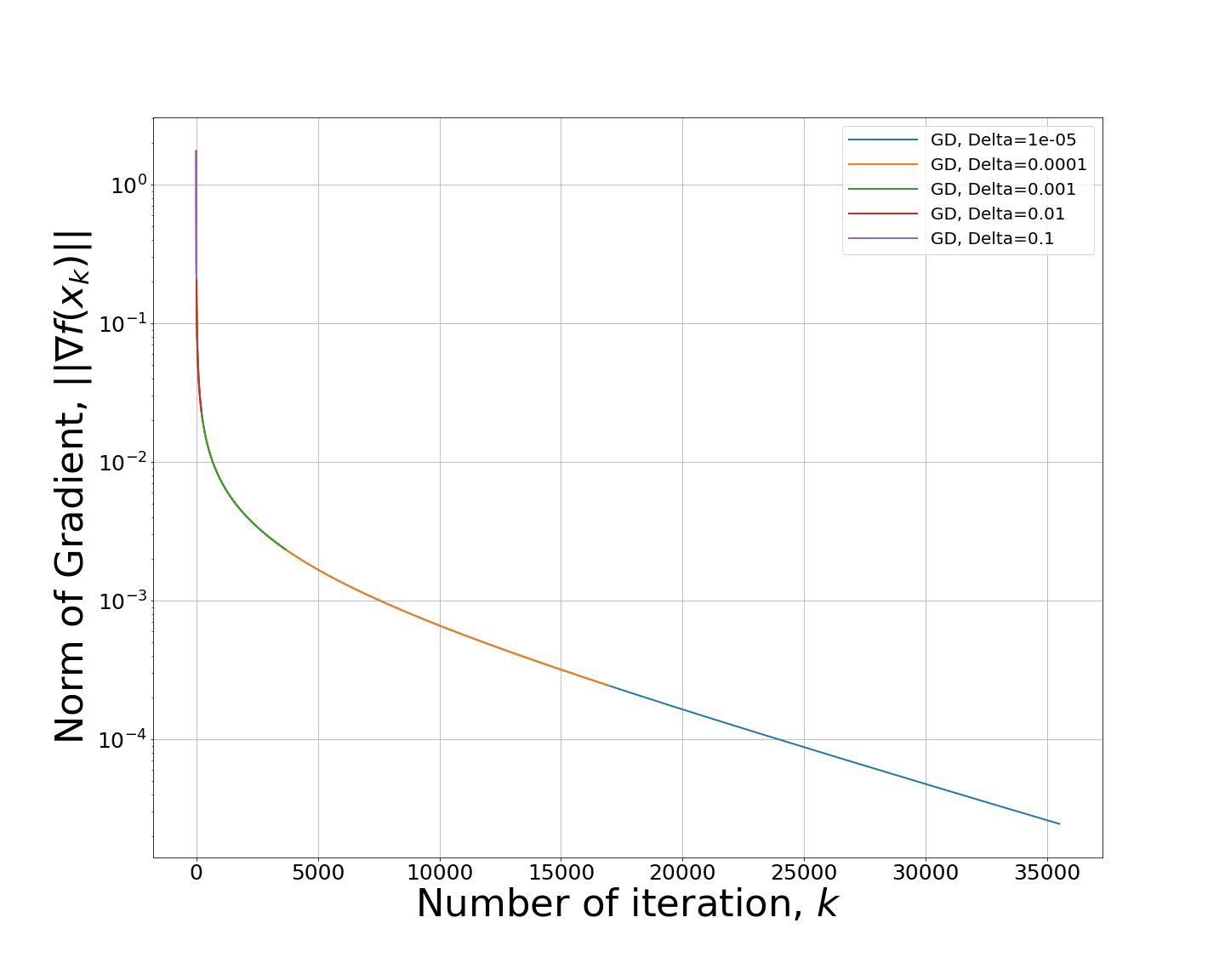}
\caption{The convergence rate of the gradient method with respect to the norm of the gradient for different values of the inexactness $\Delta$ for the problem of  logistic regression minimization with the use of  stopping criterion \eqref{f9}.} \label{fig:grad}
\end{figure}

At the same time, in Fig. \ref{fig:grad} there was shown the plot of the convergence of the gradient method with the use of stopping criterion \eqref{f9}. As we can see, the method stops when it reaches the accuracy $\|\widetilde{\nabla}f(x_k)\| \sim \Delta$. We also note, in this example, that the trajectories of the methods practically coincide until the corresponding accuracy is achieved.

\begin{figure}[ht]
\vspace{-4ex} \centering \subfigure[]{
\includegraphics[width=0.45\linewidth]{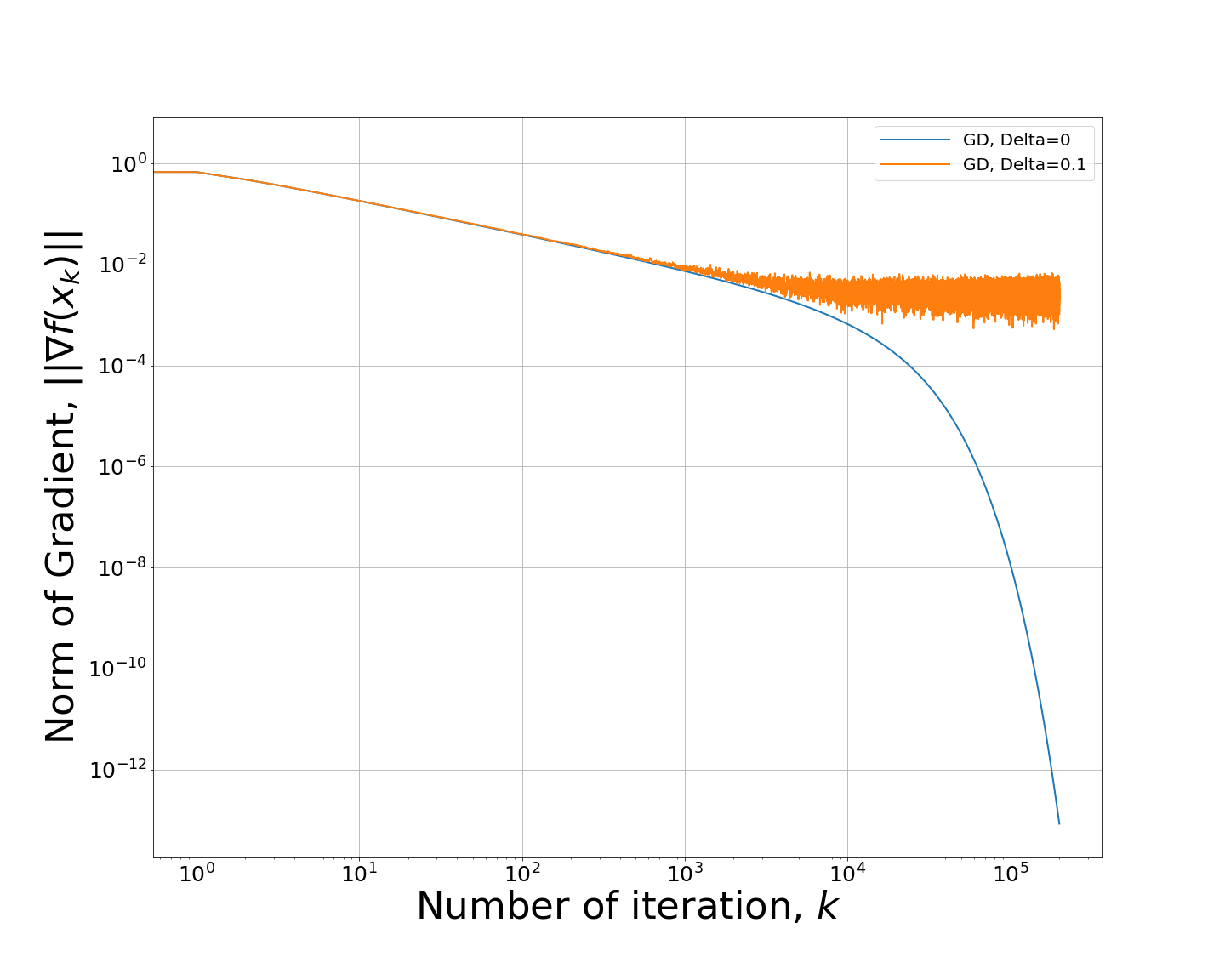}  \label{fig:grad_delta} }
\hspace{4ex}
\subfigure[]{
\includegraphics[width=0.45\linewidth]{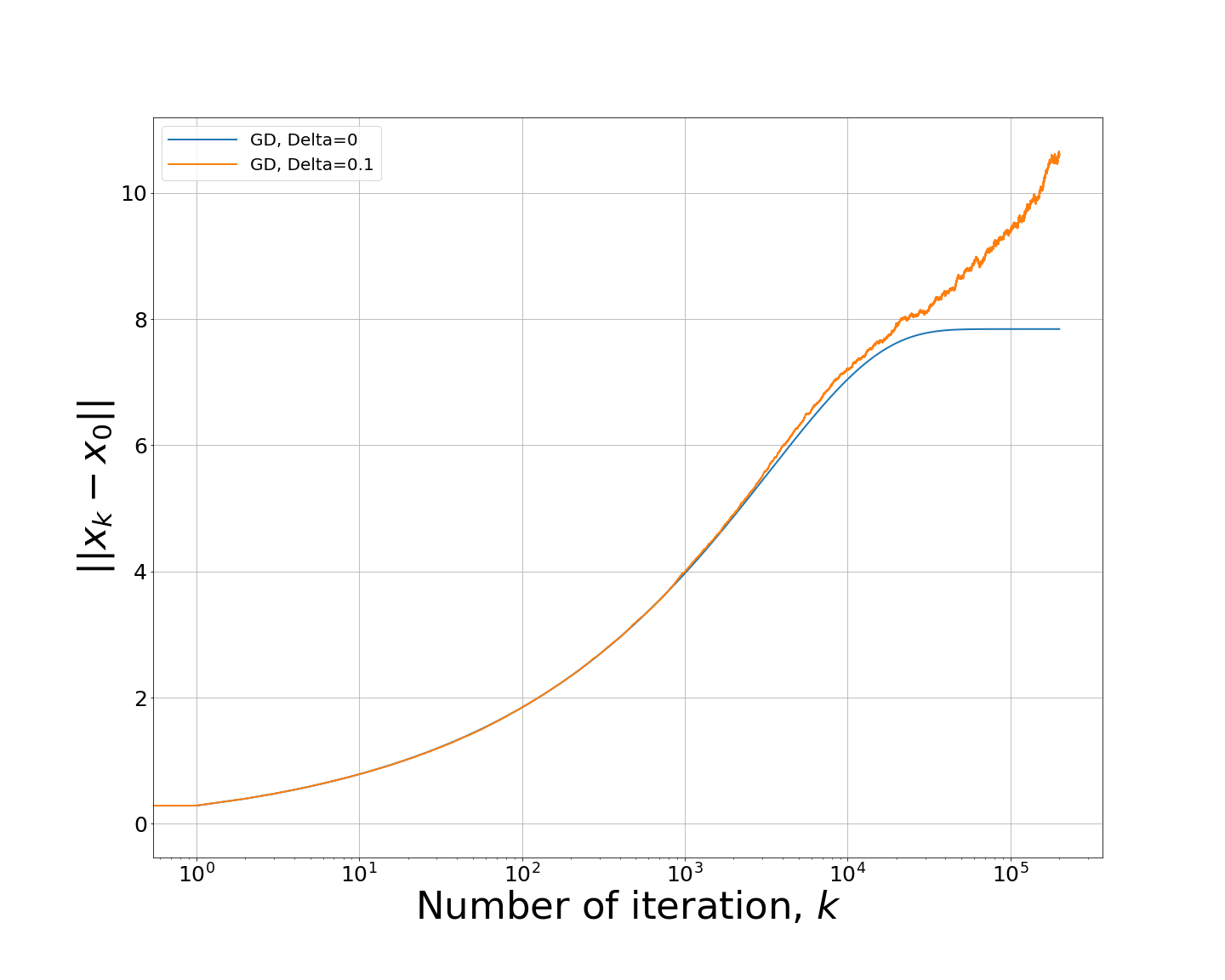}  \label{fig:norm_x_delta} }
\caption{The results of the gradient method with respect to the norm of the gradient without using the stopping criterion for $\Delta=0.1$ for the problem of logistic regression minimization for the inexactness $\Delta v(x)$. \subref{fig:grad_delta} The convergence rate with respect to the norm of the gradient; \subref{fig:norm_x_delta} the distance from the starting point to $x_k$.} \label{fig:dtype1}
\end{figure}

In Fig. \ref{fig:grad_delta} the plot of the convergence of the gradient method with the noise level $\Delta=0.1$ without using the stopping is shown. We can see that, after reaching $\|\widetilde{\nabla} f(x_k)\|\sim \Delta=0.1$, gradient method \eqref{f5} slows down significantly compared to the  noise-free method \eqref{eqf_3}. Also note that if the gradient method works without using the stopping criterion, the distance from the starting point grows uncontrollably. Moreover, this growth exceeds the distance increase for the method without noise.

\begin{figure}[ht]
\centering
\includegraphics[width=0.5\linewidth]{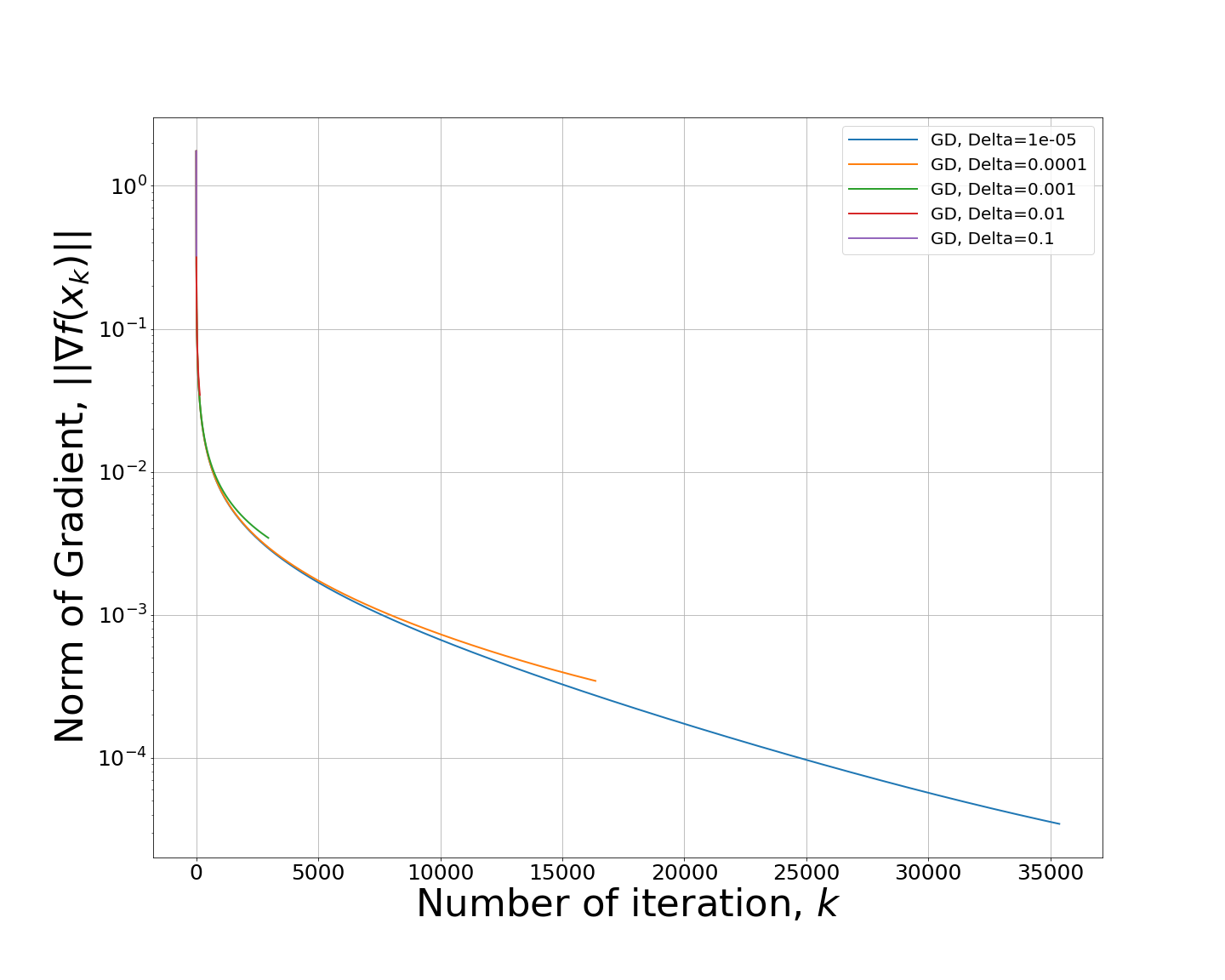}
\caption{The convergence rate of the gradient method with respect to the norm of the gradient for different values of the inexactness $\Delta$ for the problem of logistic regression minimization using stopping criterion \eqref{f9} for $v=-\frac{\Delta}{\| \nabla f(x)\|}\nabla f(x)$.} \label{fig:grad3}
\end{figure}

\begin{figure}[ht!]
\vspace{-4ex} \centering \subfigure[]{
\includegraphics[width=0.45\linewidth]{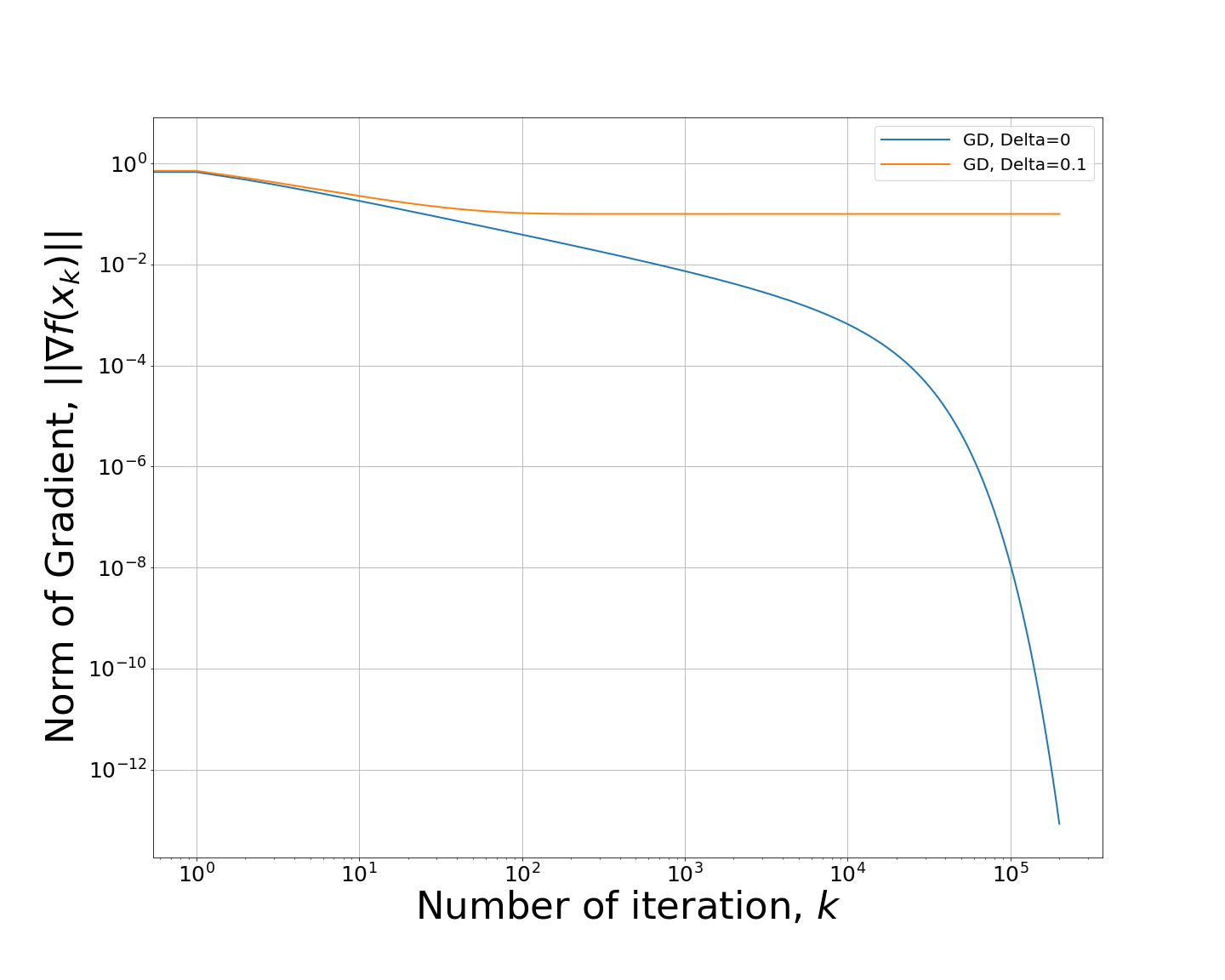}  \label{fig:grad3_delta} }
\hspace{4ex}
\subfigure[]{
\includegraphics[width=0.45\linewidth]{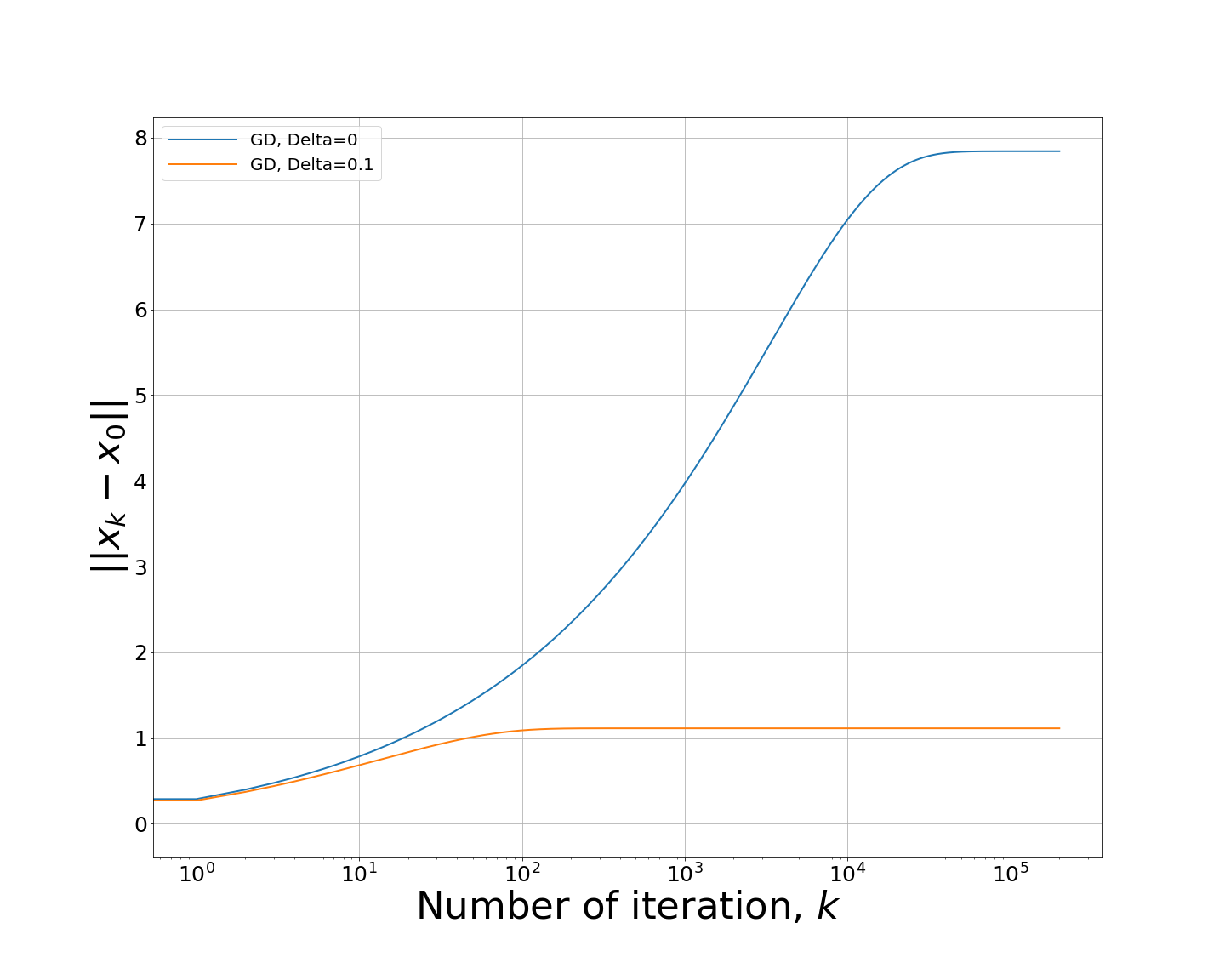}  \label{fig:norm_x3_delta} }
\caption{The results of the gradient method with respect to the norm of the gradient without using the stopping criterion for $\Delta=0.1$ for the problem of logistic regression minimization for $v=-\frac{\Delta}{\| \nabla f(x)\|}\nabla f(x)$. \subref{fig:grad3_delta} The convergence rate with respect to the norm of the gradient; \subref{fig:norm_x3_delta} the distance from the starting point to $x_k$.} \label{fig:dtype3}
\end{figure}

Next, we consider other types of problems in which the additive inexactness of the gradient proposed in section \ref{sec:exp_qp} occurs. The results for the inexact case \textbf{Antigradient} are shown in Fig. \ref{fig:grad3} and \ref{fig:dtype3}. In  Fig. \ref{fig:grad3}, we can see that the trajectories of the method begin to noticeably differ from the previous case. The inexact gradient method converges more slowly as the value of the parameter $\Delta$ increases.

However, on the graph \ref{fig:norm_x3_delta}, we can see that when methods use such an inexact gradient, the error does not accumulate. Indeed, the method stabilizes near a point such that $\|\nabla f(x)\|\approx 0.1 = \Delta$, as it can be seen from Fig. \ref{fig:grad3_delta}. Moreover, the method stops at the distance of $\|x_k-x_0\|\approx 1$ and does not move further.

Thus, in the case of $v=-\frac{\Delta}{\| \nabla f(x)\|}\nabla f(x)$, we can see that there is no problem of too large growth of the distance from the starting point to the resulting one.

The situation is essentially different for the inexactness in form of constant vector $v$, which is the same at all iterations. From Fig. \ref{fig:grad4} it can be seen that the trajectories of the method are also not the same. Moreover, adding inexactness slows down the convergence somewhat. On the other hand, the trajectories have become more similar compared to the case of the inexactness directed along the antigradient.

\subsection{Some experiments with the Rosenbrock function}
\label{sec:app_exp_rosenbrock}

In this subsection, we investigate the behavior of the proposed adaptive Algorithm \ref{adapt_gd} for  the well-known two-dimensional Rosenbrock function
$$f(x_1, x_2) = 100(x_2-(x_1)^2)^2 + (x_1-1)^2.$$

This function is not convex, and it satisfies the Lipschitz condition for the gradient only locally. Indeed, if we consider the line $x_2=0$, then we get  $f(x_1,0)=100x_1^4+(x_1-1)^2$. The gradient of this function does not satisfy the Lipschitz condition. On the other hand, the Rosenbrock function satisfies locally the PL-condition. Indeed, let us consider the system of nonlinear equations $g_1(x_1,x_2)=10(-x_2+(x_1)^2), g_2(x_1,x_2)=x_1-1$. The Jacobian of this system is
$$
    J=\begin{bmatrix}20x_1& -10\\ 1 & 0\end{bmatrix},
$$
and consequently,
$$
    JJ^\top =\begin{bmatrix}400x_1^2+100& 20x_1\\ 20x_1 & 1\end{bmatrix}\succ 0.
$$
Thus, for any compact set there exists some constant $\mu,$ such that $JJ^\top \succeq \mu I$. Then, according to the  results given in the introduction \cite{Nesterov2006}, the Rosenbrock function satisfies the PL-condition with the constant $\mu$ on the corresponding compact set.

In the conducted experiments, we vary the value of the parameter $\Delta$ and take $\delta=\Delta^2.$ The starting point for all parameters is $x_1=1, x_2 = 2$. The distance from the initial point to the optimal point $\mathbf{1} = (1,1)$ is $1$. In Table \ref{tab:rosenbrock2}, we show the results for different types of noise. The vector $v=(1,0)^\top$ was taken as a constant noise. In this experiment (and also in the next subsection), we will use stopping criterion \eqref{stop_cond_adaptGD}.

\begin{table}[ht]
    \centering
    \begin{tabular}{|c|c|c|c|c|c|c|c|c|}
            \hline
      Inexactness& $\Delta$ & Iters & Time, ms & $\|x_N-x_0\|$ & $\frac{\|\nabla f(x_N)\|}{\Delta}$ & $f(x_N)-f_*$\\
      \hline
Random & \begin{tabular}{@{}c@{}} $10^{-4}$ \\ $10^{-3}$ \\ $10^{-2}$ \end{tabular}&\begin{tabular}{@{}c@{}} $7266$ \\ $5412$ \\ $3690$ \end{tabular}&\begin{tabular}{@{}c@{}} $2273.26$ \\ $2594.61$ \\ $3163.32$ \end{tabular}&\begin{tabular}{@{}c@{}} $0.999$ \\ $0.994$ \\ $0.940$ \end{tabular}&\begin{tabular}{@{}c@{}} $2.70$ \\ $2.90$ \\ $2.61$ \end{tabular}&\begin{tabular}{@{}c@{}} $0.89 \cdot 10^{-7}$ \\ $1.00 \cdot 10^{-5}$ \\ $0.89 \cdot 10^{-3}$ \end{tabular}\\
\hline
Antigradient & \begin{tabular}{@{}c@{}} $10^{-4}$ \\ $10^{-3}$ \\ $10^{-2}$ \end{tabular}&\begin{tabular}{@{}c@{}} $7188$ \\ $5493$ \\ $3536$ \end{tabular}&\begin{tabular}{@{}c@{}} $2615.61$ \\ $2490.56$ \\ $3031.05$ \end{tabular}&\begin{tabular}{@{}c@{}} $0.999$ \\ $0.993$ \\ $0.931$ \end{tabular}&\begin{tabular}{@{}c@{}} $2.99$ \\ $3.00$ \\ $2.99$ \end{tabular}&\begin{tabular}{@{}c@{}} $0.11 \cdot 10^{-6}$ \\ $0.11 \cdot 10^{-4}$ \\ $0.12 \cdot 10^{-2}$ \end{tabular}\\
\hline
Constant & \begin{tabular}{@{}c@{}} $10^{-4}$ \\ $10^{-3}$ \\ $10^{-2}$ \end{tabular}&\begin{tabular}{@{}c@{}} $7491$ \\ $5697$ \\ $3965$ \end{tabular}&\begin{tabular}{@{}c@{}} $2301.32$ \\ $2490.32$ \\ $3485.89$ \end{tabular}&\begin{tabular}{@{}c@{}} $1.000$ \\ $0.997$ \\ $0.965$ \end{tabular}&\begin{tabular}{@{}c@{}} $1.54$ \\ $1.87$ \\ $1.93$ \end{tabular}&\begin{tabular}{@{}c@{}} $0.27 \cdot 10^{-7}$ \\ $0.24 \cdot 10^{-5}$ \\ $0.30 \cdot 10^{-3}$ \end{tabular}\\
\hline
\end{tabular}
    \caption{The results of the adaptive gradient descent for the 2D Rosenbrock function using   stopping criterion \eqref{stop_cond_adaptGD}.}
    \label{tab:rosenbrock2}
\end{table}

As previously, from the results presented  in Table \ref{tab:rosenbrock2}, we can see that the number of required iterations increases with decreasing $\Delta$ (which also tightens the stopping condition). Moreover, it increases logarithmically, which coincides with the results of Theorem \ref{theorem:adaptL_inexactf}. We can also note that the resulting distance from the starting point $x_0$ to the last point does not exceed the distance from the starting  point $x_0$  to the nearest optimal one $x_* = (1,1)$ everywhere. In addition, for all considered types of the gradient error (noise), a comparable convergence rate is observed according to the number of iterations until stopping criterion \eqref{stop_cond_adaptGD} is satisfied, and to the running time for the corresponding values of $\Delta$.

\begin{figure}[ht]
\centering
\includegraphics[width=0.5\linewidth]{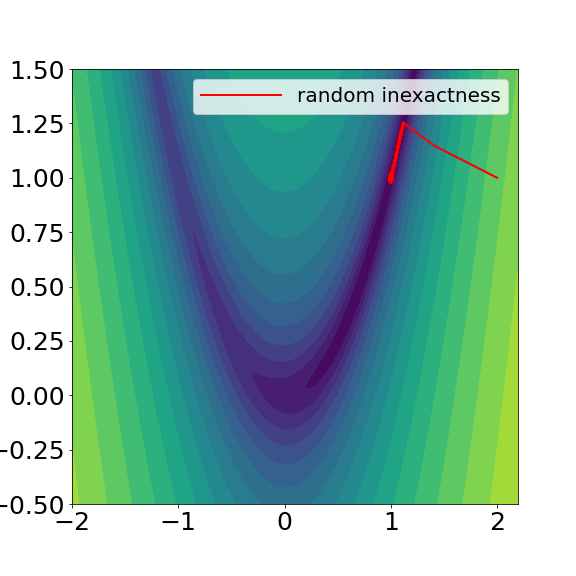}
\caption{The trajectory of the gradient method on 2D Rosenbrock function with Random Inexactness in gradient} \label{fig:rosen_traj}
\end{figure}

Note that in this example, the gradient method without a stopping criterion approaches some level for the function and next iterations are meaningless. So, in Fig.~\ref{fig:rosen_traj} we can see that the method stopped to improve the function value after some iterations.

\subsection{Some experiments with the Nesterov-Skokov function}
\label{sec:app_exp_ns}

Let us consider a system of nonlinear equations $g(x)=0,$ where $g_1=\frac{1}{2}(x_1-1), g_i=x_{i}-2x_{i-1}^2+1, i =\overline{2,n}.$ The problem of solving this system is equivalent to minimizing the following Nesterov-Skokov function (see \cite{NestSkok}):
\begin{equation}
    \label{NSfunction_app}
    f(x)=\frac{1}{4}\left(1-x_1\right)^2+\sum\limits_{i=1}^{n-1} \left(x_{i+1}-2x_i^2+1\right)^2.
\end{equation}
This function is analogous to the Rosenbrock function. It is also non-convex and satisfies the Lipschitz gradient condition only locally. Also,  function \eqref{NSfunction_app} has a global minimum at the point $(1,1\dots 1,1)^\top$ and an optimal value $f^*=0$.

Let $J$ be the Jacobian of the function $g$. Then note that $JJ^\top$ is a tridiagonal matrix, and one can easily verify that all its minors are positive for any $x$ (see the proof in appendix \ref{app:NS_function}). Whence it follows that for any compact set there exists some constant $c$ such that $JJ^\top \succeq cI$. Thus, this function locally satisfies the PL-condition.

As it was seen from the results of the previous experiments, our proposed stopping criterion \eqref{stop_cond_adaptGD} of Algorithm \ref{adapt_gd} can work equally well for all considered types of noise in the gradient. In the current experiments, for the Nesterov-Skokov function, we used the random noise of the gradient which is uniformly distributed  on the sphere. For the experiments, the starting point is $(-1,1,\dots 1, 1)^\top$ and therefore $\|x_0-x_*\|=2$. We vary the value of the inexactness $\Delta$ and the dimension of the problem $n$.

\begin{table}[ht]
    \centering
    \begin{tabular}{|c|c|c|c|c|c|c|c|c|}
      \hline
            \hline
      $n$ & $\Delta$ & Iters & Time, ms & $\|x_N-x_0\|$ & $\frac{\|\nabla f(x_N)\|}{\Delta}$ & $f(x_N)-f^*$\\
      \hline
3 & \begin{tabular}{@{}c@{}} $10^{-4}$ \\ $10^{-3}$ \\ $10^{-2}$ \end{tabular}&\begin{tabular}{@{}c@{}} $14097$ \\ $2477$ \\ $606$ \end{tabular}&\begin{tabular}{@{}c@{}} $230.58$ \\ $247.64$ \\ $383.63$ \end{tabular}&\begin{tabular}{@{}c@{}} $1.996$ \\ $2.155$ \\ $2.650$ \end{tabular}&\begin{tabular}{@{}c@{}} $2.86$ \\ $2.93$ \\ $2.19$ \end{tabular}&\begin{tabular}{@{}c@{}} $0.20 \cdot 10^{-4}$ \\ $0.11 \cdot 10^{-2}$ \\ $0.87 \cdot 10^{-2}$ \end{tabular}\\
\hline
5 & \begin{tabular}{@{}c@{}} $10^{-4}$ \\ $10^{-3}$ \\ $10^{-2}$ \end{tabular}&\begin{tabular}{@{}c@{}} $73028$ \\ $15765$ \\ $6$ \end{tabular}&\begin{tabular}{@{}c@{}} $275.03$ \\ $292.39$ \\ $200.87$ \end{tabular}&\begin{tabular}{@{}c@{}} $2.930$ \\ $3.312$ \\ $0.036$ \end{tabular}&\begin{tabular}{@{}c@{}} $2.93$ \\ $2.65$ \\ $1.45$ \end{tabular}&\begin{tabular}{@{}c@{}} $0.30 \cdot 10^{-3}$ \\ $0.49 \cdot 10^{-2}$ \\ $0.98$ \end{tabular}\\
\hline
7 & \begin{tabular}{@{}c@{}} $10^{-4}$ \\ $10^{-3}$ \\ $10^{-2}$ \end{tabular}&\begin{tabular}{@{}c@{}} $2898$ \\ $103$ \\ $17$ \end{tabular}&\begin{tabular}{@{}c@{}} $316.51$ \\ $164.23$ \\ $104.77$ \end{tabular}&\begin{tabular}{@{}c@{}} $0.049$ \\ $0.036$ \\ $0.036$ \end{tabular}&\begin{tabular}{@{}c@{}} $2.69$ \\ $2.07$ \\ $1.42$ \end{tabular}&\begin{tabular}{@{}c@{}} $0.98$ \\ $0.98$ \\ $0.98$ \end{tabular}\\
\hline
\end{tabular}
    \caption{The results of the adaptive gradient descent for the Nesterov-Skokov function with the use of stopping criterion \eqref{stop_cond_adaptGD}.}
    \label{tab:NestSkok}
\end{table}

Table \ref{tab:NestSkok} shows the results of the adaptive gradient method \ref{adapt_gd} for the Nesterov-Skokov function \eqref{NSfunction_app}. Firstly, we see that as the dimension of $n$ increases, the difference between the required time to solve the problem for different $\Delta$ grows significantly. Secondly, for different $n$ with the same $\Delta$, the method converges to a solution with significantly different accuracy. So for $\Delta=10^{-4}$ the accuracy for $n=7$ and $n=3$ differs by more than 100 times. This is explained by the decrease in the constant $\mu$ in the PL-condition as $n$ grows. We can also note that $\|x_N-x_0\|$ exceeds $\|x_0-x_*\|$ by at most 2 times. Moreover, significant upward deviations are observed for the cases when numerous iterations are made ($n=5$ and $\Delta=10^{-4}, 10^{-3}$). It can also be noted that even for sufficiently small values of the norm of the gradient, the accuracy by the function turns out to be quite low (which is typical for the Nesterov-Skokov function). Thus, for $n=7$ and $\Delta=10^{-4}$ we get a point with $\|\nabla f(\widehat{x})\|\approx 10^{-4},$ but $f(\widehat{x})-f^*\approx 0.98$.

\begin{figure}[ht!]
\vspace{-4ex} \centering \subfigure[]{
\includegraphics[width=0.45\linewidth]{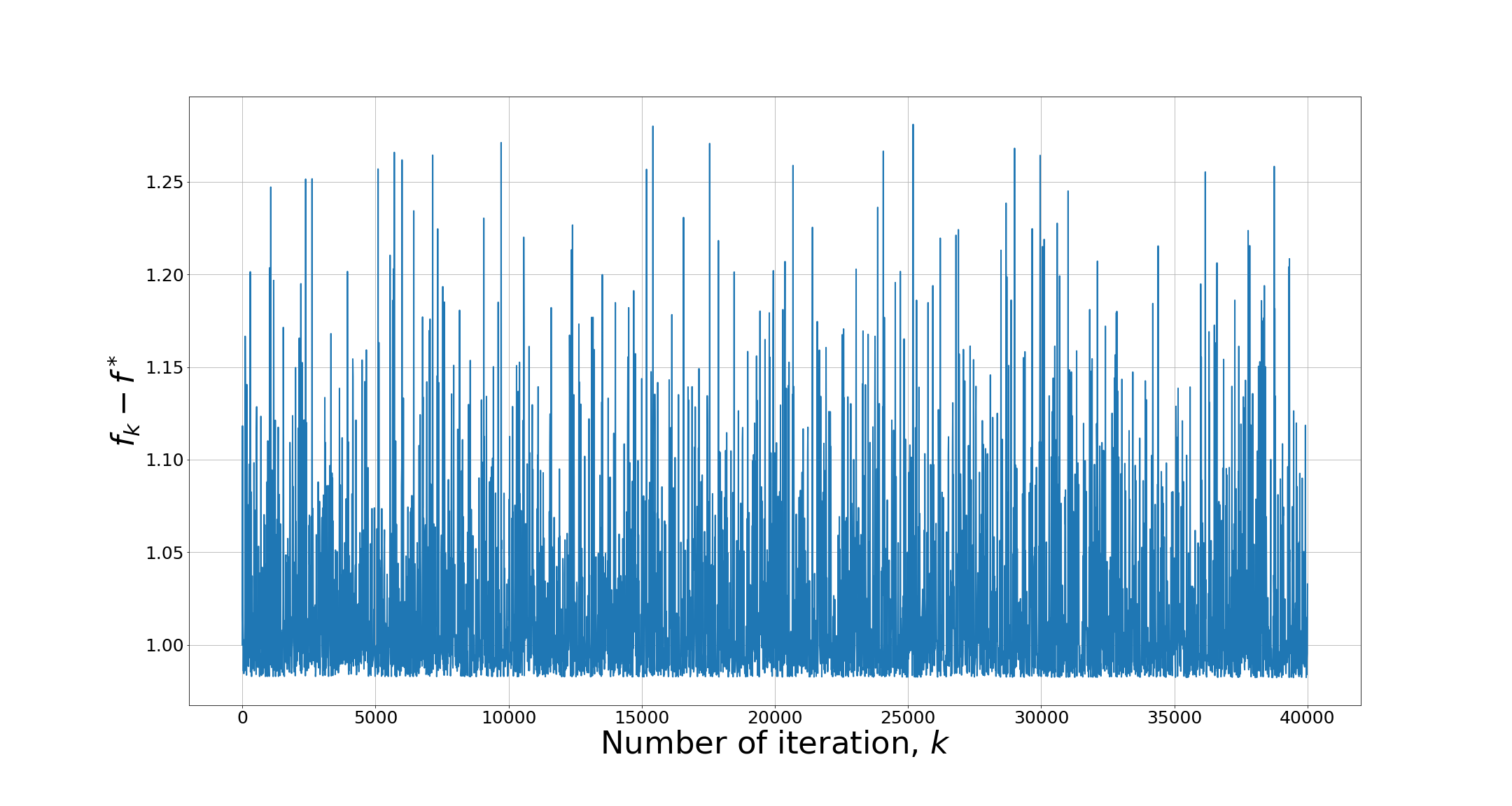}  \label{fig:NS_traj_value} }
\hspace{4ex}
\subfigure[]{
\includegraphics[width=0.45\linewidth]{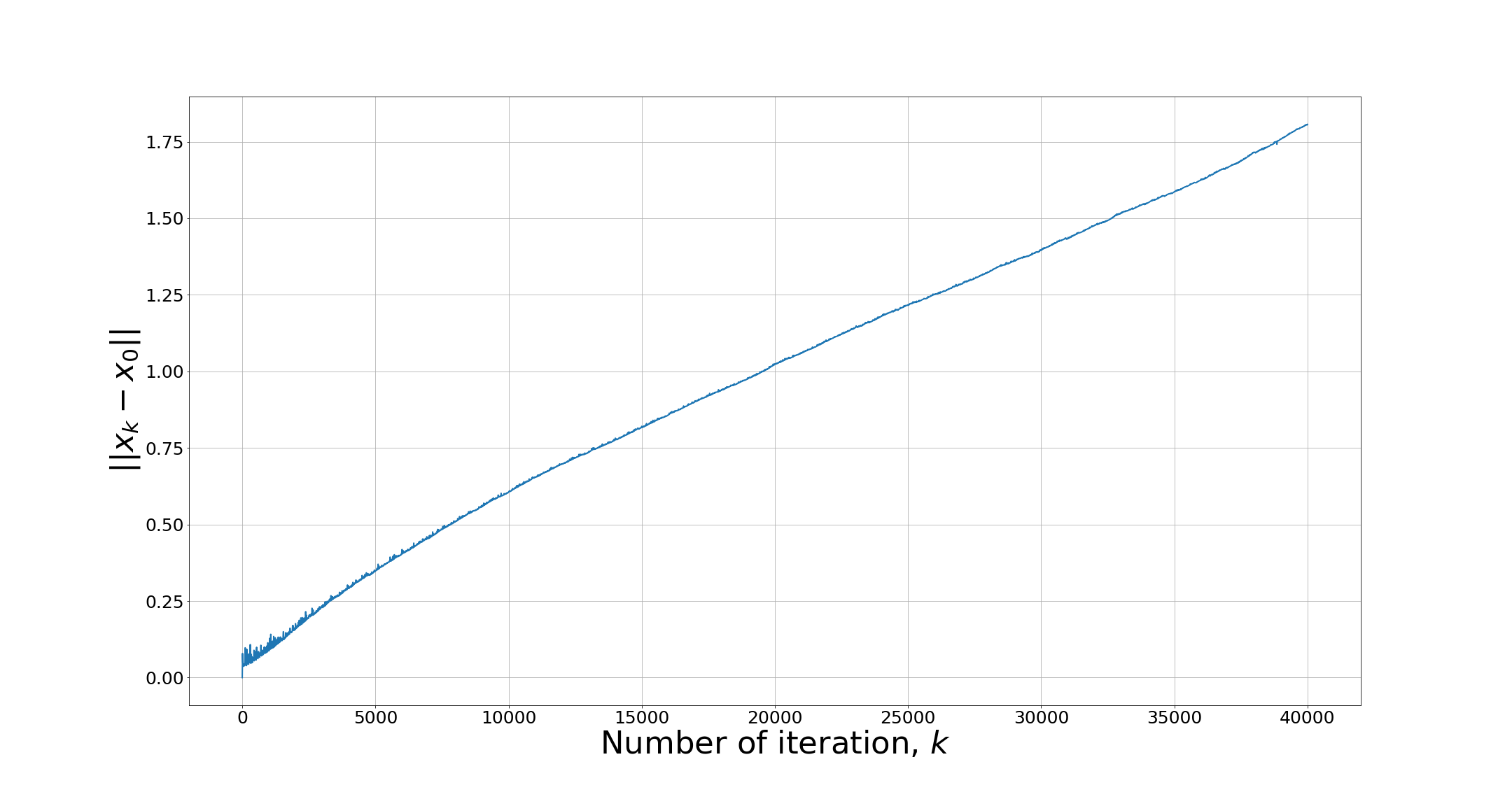}  \label{fig:NS_traj_dist} }
\caption{The convergence of the gradient method for the Nesterov-Skokov function in the 7 dimensional space with the Antigradient inexactness in the gradient: \subref{fig:NS_traj_value} the function value $f(x_k)-f^*$; \subref{fig:NS_traj_dist} the distance $\|x_k-x_0\|$.} \label{fig:NS_traj}
\end{figure}

In Fig. \ref{fig:NS_traj} we can see the convergence rate of the gradient method with the Antigradient inexactness. First, we can see in Fig. \ref{fig:NS_traj_value} that the method approaches level 0.98 quite quickly and does not improve after that. On the other hand, from  Fig. \ref{fig:NS_traj_dist} one can see that the method without the stopping criterion moves away from the starting point quite far.

\section{Proof of Proposition \ref{prop:logreg}}
\label{proof_prop_logreg}

Further, $x_* \in X_*$. Condition 1 in  Proposition \ref{prop:logreg} holds to any logistic regression function. Let us estimate the Lipschitz constant  $L$ of the gradient of  function \eqref{log_regr}. It is known that 
$
L = \max_{x\in\mathbb{R}^n} \lambda_{\max}\left( \nabla^2 f(x)\right).
$
On the other hand, $\nabla^2 f(x)= \nabla^2 g(Ax) = A^\top H_g(z)\Big|_{z=Ax} A,$
where $g(z)=\frac{1}{m}\sum\limits_{i=1}^m\log \left(1 + \exp\left(z_i\right)\right)$, and $A= [-y_1 w_1\dots -y_m w_m]^\top \in \mathbb{R}^{m \times n}$ and the matrix $H_g(z)$ is the Hessian of the function $g$ at the point $z$. Note that $H_g(z)$ is a diagonal matrix with entries $\frac{1}{m}\frac{e^{z_i}}{(e^{z_i}+1)^2}\leq \frac{ 1}{4m}$. Thus, we have the following estimate from above:
$$
L\leq \|A\|_2^2 \max_{z}\|H_g(z)\|_2=\frac{\|A\|_2^2}{4m}.
$$
So,  statement 2 in Proposition \ref{prop:logreg} holds.

Further, let us introduce new notations. Let $E_1$ be a subspace given by the basis $W_B$ and $E_2$ be a subspace orthogonal to $E_1$. Note that if $k<n$, then the dimension of $E_2$ is at least 1. Then if there exists a minimum point $x_*$, then at any point from the set $x_*+E_2\subseteq X^*$, the objective function takes the minimal value. Therefore, the set of the solutions is unbounded.

Now let us prove statement 3 in  Proposition \ref{prop:logreg}. Note that the created matrix $W$ has a rank $k\leq n$. Accordingly, all vectors $w_i$ belong to the $k$-dimensional subspace $E_1$, given by the basis $W_B$. Also, for any vector $\widetilde{x}\in E_2$ from the subspace orthogonal to $E_1$ and for any vector $x\in\mathbb{R}^n$, it is true that $f(x+\widetilde{x} )=f(x)$. Thus, $f(E_1)=f(\mathbb{R}^n)$. Note that function \eqref{log_regr} is bounded from below by 0, and hence $f^*=\inf_{x} f(x)\geq 0>-\infty$. Thus, we consider a sequence $\{x_j\}_j\in E_1$ such that $f(x_j)\rightarrow f^*$. Let us transform the sum, taking into account that the first $2k$ vectors $w_j$ are rows of the matrices $W_B$ and $-W_B$
\begin{align*}
\sum\limits_{i=1}^m\log \left(1 + \exp\left(-y_i \langle w_i, x_j\rangle\right)\right)&  \geq \sum\limits_{i=1}^{2k}\log \left(1 + \exp\left(-y_i \langle w_i, x_j\rangle\right)\right) \\
&=\sum\limits_{i=1}^{k}\log\left( \left(1 + e^{-y_i \langle w_i, x_j\rangle}\right)\left(1 + e^{y_i \langle w_i, x_j\rangle}\right)\right)\\
&=\sum\limits_{i=1}^{k}\log\left(2+2\text{ch}\left(y_i \langle w_i, x_j\rangle\right)\right).
\end{align*}

From the fact that $f^*$ is  finite  and from the constructed lower bound, it follows that $|\langle a_i, x_j\rangle|\leq C, \forall j$ for some constant $C > 0$, i.e. $\|W_B x_j\|\leq k C, \; \forall j$.  On the subspace $E_1$, the matrix $W_B$ defines an invertible operator. Hence, $\|\cdot\|_{W_B^\top W_B}$ is a norm on the subspace $E_1$. Therefore, in view of the equivalence of norms in a finite-dimensional space, we have that $\|x_j\| \leq C_1 \forall j$ for some constant $C_1>0$ depending only on the constants $C,k$ and the parameters of the matrix $W_B$. Thus, any sequence of elements of the space in which the sequence of values of the function converges to $f^*$ is bounded.  This means that from this bounded sequence, we can extract a convergent subsequence $\{x_{j_l}\}_l$. The limit of this subsequence is the desired point $x_*$, which is a finite vector. So, statement 3 in Proposition \ref{prop:logreg} holds.

As mentioned before, for any vector $\widetilde{x}\in E_2$ and for any vector $x\in\mathbb{R}^n$, the equality $f(x+\widetilde{x} )=f(x)$ holds. By constructing $W_B$, the dimension of $E_2$ is at least 1. So, from this and statement 3, we have that statement 4 in Proposition \ref{prop:logreg} holds.

\section{The Nesterov-Skokov function }
\label{app:NS_function}

Let us consider the following known Nesterov-Skokov function \cite{NestSkok}
\begin{equation}\label{NSfunctionapp}
    f(x)=\frac{1}{4}\left(1-x_1\right)^2+\sum\limits_{i=1}^{n-1} \left(x_{i+1}-2x_i^2+1\right)^2=\sum\limits_{i=1}^n g_i^2(x),
\end{equation}
where $g_1=\frac{1}{2}(x_1-1), g_i=x_{i}-2x_{i-1}^2+1, i=\overline{2,n}$. Then the Jacobian of the system $g(x)=0$ is a two-diagonal matrix. On the main diagonal $J_{11}=\frac{1}{2}$ and $J_{ii}=1$, on the side diagonal $J_{i,i-1}=-2x_{i-1}$. Then the matrix $JJ^\top$ is a tridiagonal symmetric matrix of the form
$$
\begin{bmatrix}
\frac{1}{4} & -2x_1 & 0 & 0 & \dots & 0 & 0 \\
-2x_1 & 16x_1^2+1 & -4x_2 & 0 & \dots & 0 & 0 \\
0 & -4x_2 & 16x_2^2+1 & -4x_3 & \dots & 0 & 0 \\
\hdotsfor{7}\\
0 & 0 & 0 & 0& \dots & 16x_{n-2}^2+1 & -4x_{n-1} \\
0 & 0 & 0 & 0& \dots & -4x_{n-1} & 16x_{n-1}^2+1
\end{bmatrix}.
$$
Then the first principal minor is $f_1=\frac{1}{4}$, the second is $f_2 = \frac{1}{4}$. The recursive formula for a tridiagonal matrix is $f_k=(16x_{k-1}^ 2+1)f_{k-1} - 16x_{k-1}^2 f_{k-2}$. Then we can prove that $f_j\geq f_{j-1}$ for all $j>1$. Thus, the matrix is strictly positive definite for any $x$. Then for any compact set $W$, one can choose a constant $c$ which limits from below all the eigenvalues of the matrix $JJ^\top$, which means that $JJ^\top \succeq cI$. Then, according to the results \cite{Nesterov2006} given in the introduction, the function satisfies the PL-condition with a constant $c > 0$ on the corresponding compact set.

\end{document}